\input amstex 
\documentstyle{amsppt}
\magnification 1100
\NoBlackBoxes
\hsize= 16.85 truecm
\vsize= 21.75 truecm
\hoffset= -.15 truecm
\voffset= 1 truecm

\def\ns{\hskip -.15 truecm}
\def\numfree{Lemma 2.13 }
\def\ab{Theorem 4.1 }
\def\ablemma{Lemma 4.2 }

\def\abfive{(4.3.1) }
\def\absix{(4.3.2) }
\def\abseven{(4.3.3) }
\def\abeight{(4.3.4) }
\def\sequence{(1.1) }

\def\CM{(1.4.2) }
\def\sfreeobs{Observation 1.4.1 }
\def\abcor{Corollary 4.4 }
\def\defK{Definition 3.1 }
\def\KLaz{Lemma 3.2 }
\def\Ksequence{Lemma 3.4 }
\def\defrK{Definition 3.3 }
\def\EKoszul{Theorem 3.5 }
\def\abKempf{Remark 4.6 }

\def\preg{Theorem 1.3 }
\def\twop+four{Corollary 5.11 }
\def\twoK{Theorem 5.5 }
\def\Butres{Proposition 2.4 }
\def\abKoszul{Theorem 4.7 }
\def\abKlemma{Lemma 4.8 }
\def\preglemma{Lemma  1.4 }

\def\foursections{Example 6.5 }

\def\slopelemma{Lemma 2.10 }
\def\Splemma{Lemma 2.9 }
\def\Splemmaone{2.9.1 }
\def\Enriquesng{ Theorem 2.2 }
\def\Enriquesnpres{Theorem 2.11 }
\def\EnriquesMukaing{Corollary 2.8 }
\def\EnriquesMukainpres{Corollary 
 2.12 }
\def\Enriqueshighersyzygies{Theorem 2.14 }
\def\MP{0.1 }
\def\Enriquessame{Theorem 2.1 }
\def\abMukai{Corollary 4.5 }
\def\twoamplearefree{Lemma 2.7 }
\def\GLlemma{Theorem 1.2 }
\def\Pareschilemma{Proposition 2.5 }
\def\robs{Observation 2.3 }
\def\EnriquesMukaiNp{Corollary 2.15 }
\def\fourK{Corollary 5.6 }
\def\fiveK{Corollary 5.9 }

\def\regulargtNp{Theorem 5.12 }

\def\Bombieri{[Bo] }
\def\Butler{[Bu] }
\def\Catanese{[Ca] }
\def\Ciliberto{[Ci] }
\def\CosDol{[CD] }
\def\Busto{[FdB] }
\def\EL{[EL] }
\def\FV{[FV] }
\def\GPone{[GP1] }
\def\GPtwo{[GP2] }
\def\GPthree{[GP3] }
\def\GPfour{[GP4] }
\def\Greentwo{[G1] }
\def\Greenthree{[G2] }
\def\Greenfour{[G3] }
\def\Hommaone{[H1] }
\def\Hommatwo{[H2] }
\def\Kaw{[Ka] }
\def\Kempf{[Ke] }
\def\Lazarsfeld{[L] }
\def\Miyaoka{[Mi] }
\def\Mumford{[Mu] }
\def\Pareschione{[P1] }
\def\Pareschitwo{[P2] }
\def\PP{[PP] }
\def\Reider{[R] }
\def\Sho{[S] }
\topmatter
\title  Very ampleness and higher syzygies for algebraic surfaces
and Calabi-Yau threefolds
\endtitle
\author Francisco Javier Gallego \\ and \\ B. P. Purnaprajna
\endauthor
\endtopmatter
\document
\bigskip
\bigskip
\bigskip
{\bf Contents.}
\bigskip
\bigskip

Abstract
\vskip .2 truecm 
{\bf Part 1.}
 
 Projective normality and syzygies of algebraic surfaces
\vskip .2 truecm
\roster
\item "" Introduction
\item "1." A general result on syzygies of algebraic varieties
\item "2." Syzygies of Enriques surfaces
\item "3." Koszul rings of Enriques surfaces
\item "4." Abelian and bielliptic surfaces
\item "5." Surfaces of positive Kodaira dimension
\item "" References
\endroster
\vskip .2 truecm

{\bf Part 2.}
 
Very ampleness and higher syzygies of Calabi-Yau threefolds
\roster
\item "" Introduction
\item "1." Very ampleness and projective normality
\item "2." Normal presentation, Koszul rings and higher syzygies
\item "" Appendix: Singular Calabi-Yau threefolds
\item "" References
\endroster

\newpage
\
\vskip 3 truecm
\heading  ABSTRACT \endheading
\bigskip

This work consists of two parts.
In the first part we develop new techniques to compute Koszul
cohomology groups for several classes of varieties. As 
applications  we prove  results on
projective normality and syzygies for 
algebraic surfaces. From
more general results we obtain in particular the following:
\roster
\item"a)"  Mukai's conjecture (and stronger
variants of it) regarding projective normality and
normal presentation for surfaces with Kodaira dimension
$0$, and
 uniform bounds for higher syzygies
associated  to adjoint linear series,
\item"b)" effective bounds along the lines of Mukai's conjecture
regarding projective normality and normal presentation for surfaces
of positive Kodaira dimension, and,
\item"c)"  results  on projective normality for pluricanonical
models of surfaces of general type (recovering and strengthening
results by Ciliberto; cf.
\Ciliberto \ns) and generalizations of them to higher syzygies.
\endroster
\medskip

In the second part  we prove results on very ampleness, projective normality and higher
syzygies for both singular and smooth Calabi-Yau threefolds. The general results
that we prove are analogues of the results of St. Donat for K3 surfaces. From these
general results we obtain bounds very close to Fujita's conjecture regarding very
ampleness of powers of an ample line bundle $A$ (for instance, if $A^3 > 1$, our
bound is one short of Fujita's). We generalize our results on very ampleness and
projective normality to higher syzygies. 

\newpage

\topmatter
\title Part 1: Projective normality and syzygies of algebraic surfaces
\endtitle
\endtopmatter
\document 
\vskip .3 cm

\headline={\ifodd\pageno\rightheadline \else\leftheadline\fi}
\def\rightheadline{\tenrm\hfil \eightpoint PROJECTIVE NORMALITY
AND SYZYGIES OF ALGEBRAIC SURFACES
 \hfil\folio}
\def\leftheadline{\tenrm\folio\hfil \eightpoint F.J. GALLEGO \& B.P. PURNAPRAJNA \hfil}

\heading  Introduction\endheading
In this article we develop  
new techniques to
compute Koszul cohomology groups.
Koszul cohomology is  important because of its
relation to Hodge Theory and to the computation of
syzygies of projective varieties. In the present work we focus on
the latter application. The topic of syzygies is  interesting
because it deals with the interplay between algebra and
geometry: 
the algebra coming from the
equations defining the variety and the geometry arising
from the knowledge of what line bundles live on the variety. The
earliest result typical of this application we have in mind is the
result of Castelnuovo, who showed that a curve of degree greater
than 
$2g$  has a normal homogeneous
coordinate ring ($g$ denotes the genus of the curve). He
also proved that if the degree was greater than $2g+1$,
then the ideal of the curve was generated by quadratic
equations. This result was rediscovered later by many
people, among others Fujita, St. Donat, Mumford, Green,
etc. Recently Mark Green threw new light on this
connection between algebra and geometry by generalizing
the study of homogeneous coordinate rings and ideals to
the study of free resolutions. He linked the behavior of  
graded Betti numbers of the resolution of the
homogeneous coordinate ring to the
cohomology groups of certain vector bundles on the
variety (see \Greentwo \ns, \Greenthree and \Greenfour \ns;  for a
particularly nice introduction to the subject see also \Lazarsfeld 
 and  for the precise statement used in this article see
\GLlemma
\ns).  Green generalized  Castelnuovo's result
proving that if the degree of the curve is greater than $2g+p$, 
then the resolution is in addition linear until the $p$th stage.
This property of the resolution is the so-called property $N_p$. 
Connection between algebra and geometry is better seen in the
case of the canonical curve. Here there are classical
results by N\"other and Petri on  projective
normality and normal presentation for canonical curves.
The geometric part of the statements is
summed up in the Clifford index of the curves. Green's
conjecture for canonical curves generalizes N\"other and
Petri's results, claiming  that the shape of the free 
resolution of the canonical ring is
determined by the Clifford index of the curve (precisely, if the Clifford
index is $p+1$, then the resolution satisfies exactly the property $N_p$). 

\vskip .2 truecm
There are still many open questions regarding linear series on
curves, but for surfaces and higher dimensional varieties the
field is almost entirely open. Among the  open questions for
surfaces and higher dimensional varieties,  the conjectures 
of Fujita on very ampleness and Mukai on higher syzygies 
of surfaces have
attracted attention in recent years.  Fujita conjectured that on an
algebraic variety $X$ of dimension $n$, if $A$ is an ample line
bundle on $X$, then $K_X
\otimes A^{\otimes n+2}$ should be very ample, where $K_X$
denotes the canonical bundle on $X$.  Mukai's conjecture
says that if $S$ is a surface, $A$ is an ample line bundle on $S$,
$L$ is a line bundle on
$S$ equal to
$K_S
\otimes A^{\otimes n}$,
and  $n \geq p+4$, then $L$ satisfies property
$N_p$. This conjecture can be regarded as a two dimensional
analogue of Green's theorem for curves.  Indeed, Green's theorem
can be interpreted as follows: any line bundle
$L$ on a curve
$C$ which is at least as positive as $ K_C \otimes A^{\otimes
p+3}$ satisfies property $N_p$, where $K_C$ is the canonical
bundle of
$C$ and
$A$ is an ample line bundle on $C$. 
Fujita's conjecture has been proved for algebraic surfaces and
it follows from a remarkable result of
 Reider (cf. \Reider \ns). For higher dimensional varieties some
effective  bounds have been obtained. Even though 
the bounds are far from what has been conjectured, 
they are considered an important
step towards the goal of proving  Fujita's conjecture.
Mukai's conjecture has not yet been proved even
for $p=0$.
Some progress has been made  
by Butler for ruled varieties
(see \Butler \ns), 
Kempf for Abelian varieties (see
\Kempf \ns), and Ein and Lazarsfeld, who prove a beautiful, very
general result on  
adjoint linear series associated to 
very ample line bundles (see \EL \ns).  Y. Homma  proved Mukai's conjecture 
for the case $p=0$ for elliptic ruled surfaces (see \Hommaone and 
\Hommatwo
\ns). One of the things we do here is
to prove Mukai's conjecture in certain cases and obtain effective bounds 
towards it for all surfaces.
\vskip .2truecm
In this article we pursue a new direction to study syzygies of
algebraic surfaces. This direction can be summarized in the 
following  meta-principle:

\proclaim{\MP} If $L$ is  the product of $(p+1)$ ample and 
base-point-free line bundles
satisfying ``certain cohomological" conditions, then $L$ satisfies 
the condition $N_p$.
\endproclaim

\noindent With the meta-principle as a guiding light,
we obtain the following as  corollaries of our more general results:
\roster
\item"(1)" We prove that Mukai's conjecture regarding projective 
normality and normal presentation  is true,   lowering Mukai's
bound by one in the latter case,
for all   surfaces of Kodaira dimension $0$ and show stronger variants
of it (cf. corollaries 2.6 and 4.5).  

\item"(2)" We  obtain a uniform bound along the line of Mukai's
conjecture for higher syzygies associated  to adjoint linear series for
all   surfaces of Kodaira  dimension
$0$.

\item"(3)" We find effective bounds along the lines of Mukai's
conjecture regarding projective normality and normal presentation
for surfaces of positive Kodaira dimension.

\item"(4)" We obtain results on
projective normality, normal presentation  and higher syzygies 
for pluricanonical models
of surfaces of general type, 
 recovering and strengthening
results of  Ciliberto (cf. \Ciliberto \ns).

\item"(5)" We find effective bounds regarding projective normality
and higher syzygies for multiples of ample bundles for regular
surfaces of positive Kodaira dimension.
\endroster

\noindent Result (3) can be interpreted as a
higher syzygy  analogue along the lines of Mukai's of
the effective results of
Demailly, Ein and Lazarsfeld regarding Fujita's conjecture.
Result (5) is a higher syzygy  analogue of the results of Siu and
Fern\'andez del Busto regarding effective Matsusaka's theorems on
base-point-freeness and very ampleness. 
On the other hand, since an 
effective bound regarding Mukai's conjecture was obtained by Butler
for ruled surfaces,  results (1),  (2)  and (3), coupled with
Butler's give the best bounds so far towards Mukai's conjecture for
all surfaces.
\vskip .2 truecm
Almost all known
results on syzygies of algebraic surfaces (and several results on 
curves) fit into 
\MP \ns. For example, the normal presentation of line bundles
of  degree greater than
$2g+1$ on curves, by Castelnuovo and others (see
\text{\GPone \ns),} the result of Kempf referred to above
(see \text{\abKempf \ns),} and the result by Ein and
Lazarsfeld.  In \GPone \ns, \GPtwo and \GPthree we show the
validity of
\MP  for surfaces of Kodaira dimension
$-\infty$ and K3 surfaces. We 
show in the present article that
\MP holds for all other surfaces of Kodaira dimension $0$ and
for adjoint linear series (more general than those involved in
Mukai's conjecture) on surfaces of positive Kodaira dimension.
We summarize here the results which give evidence of the above
claims: 

For surfaces of Kodaira dimension $-\infty$, 
 the $(p+1)$-th power of an ample, free and
nonspecial line bundle satisfies property $N_p$ (\GPtwo
\ns, Theorem 2.2, see also Lemma 2.8; our result is in fact
more general as it is stated for surfaces with geometric genus
$0$).  \preg of this paper
generalizes this result and  unifies among others  
\twop+four for surfaces of general type and Corollary 1.6
for Calabi-Yau threefolds.  
Moreover, in
\GPone and
\GPtwo we prove finer versions  (\text{\GPone\ns,} Theorem 4.2 and
\GPtwo \ns, Theorem 6.1) of the meta-principle for elliptic ruled
surfaces, yielding among other things the fact that Mukai's
conjecture regarding normal presentation holds for such
surfaces. For anticanonical rational surfaces we also prove finer
versions of \MP \ns, in a modified version of \GPthree \ns.

For surfaces with Kodaira dimension $0$ we
show precisely the following:
\vskip .15 truecm
\proclaim{Theorem 0.2} Let $S$ be a minimal surface with Kodaira
dimension
$0$ and let $B_1, \dots, B_n$ be  numerically equivalent,
ample and base-point-free line bundles. Assume that the sectional 
genus  of the $B_i$ is greater than or equal to $2$ if $X$ is an
Enriques surface and that the $B_i$ are non-hyperelliptic with
sectional genus greater than or equal to $4$ if
$X$ is a K3 surface. Then $B_1
\otimes \cdots \otimes B_n$ satisfies $N_p$ for all
$n \geq p+1$ and $p \geq 1$.
\endproclaim
\noindent The proof of Theorem 0.2  can be found for Enriques surfaces 
in Section 2, for Abelian and bielliptic surfaces in Section 4 
and for K3 surfaces in \GPthree \ns. In the case of K3 surfaces 
we  prove a stronger version of \MP imposing extra conditions on
$B_i$ (see
\text{\GPthree \ns).} 
As a consequence of Theorem 0.2 we obtain the following:
\proclaim{Theorem 0.3}
Let $S$ be a minimal surface with Kodaira
dimension
$0$, let $B$ be an ample and base-point-free line bundle,
and let $A$ be an ample line bundle. If $n
\geq p+1$ and $p \geq 1$, then the bundle $K_S
\otimes B^{\otimes n}$ satisfies property $N_p$ and  if $m \geq
2p+2$ and $p \geq 1$, then the bundle
$K_S \otimes A^{\otimes m}$ satisfies property $N_p$. 
\endproclaim

Theorem 0.3 recovers Kempf's result for Abelian surfaces and
implies the already mentioned result (1) regarding Mukai's
conjecture.  

For surfaces of positive Kodaira dimension we prove
results in the spirit of \MP for adjoint linear series and for 
powers of ample 
and  base-point-free line bundles (see Theorems 5.1,
5.8 and 5.14). We apply these results to obtain
the above mentioned results
(3) and (5) regarding
 effective bounds for projective normality, normal presentation,
and property $N_p$, and (4)  on
 pluricanonical models
of surfaces of general type (see \regulargtNp
and
Theorem 5.16). 
Our results on projective normality, normal presentation and higher 
syzygies of pluricanonical models recover and strengthen
results of  Ciliberto on projective normality. In particular, we
show the following, which is a question  posed by Bombieri (in
[Bo]):

\vskip .15 truecm 
\noindent {\it Let $X$ be a
surface of general type such that $p_g \geq 2$ or $K_X^2
\geq 5$. If $n \geq 5$, then the image of $X$
by $|K_X^{\otimes n}|$ is projectively normal}.
\vskip .15 truecm
\noindent Moreover we improve results of [Ci] in the case of 
regular surfaces (\fourK \ns). 

\vskip .2 truecm

We apply the techniques developed in this article to study syzygies of higher
dimensional varieties. We show results in the spirit of \MP  for
Fano varieties in \GPthree and for Calabi-Yau threefolds.
In \GPfour we prove optimal results
on very ampleness, projective normality and higher syzygies for Calabi-Yau
threefolds. These results are similar in spirit
to the well known results of St. Donat for K3 surfaces
and Lefschetz for Abelian varieties. 
\vskip .2 truecm
Another very interesting problem in this area is the
relation  between 
normal presentation and the Koszul property of coordinate rings.  
We show that whenever a line bundle  on the variety under
consideration (in this article) is normally presented then it
embeds the variety with a Koszul homogeneous coordinate ring. This
gives further evidence to the following  (to paraphrase Arnold):
Any homogeneous coordinate ring which has a serious reason for
being quadratically presented is Koszul. In  Section 3 we develop
the necessary tools to tackle  this problem and restrict
ourselves to Enriques surfaces. In the subsequent sections we
apply these tools to prove the result for other surfaces.

\vskip .2 truecm
 A basic obstacle one encounters in the  kind of problems we have
been talking about in the previous paragraphs is the scarcity of
techniques to compute Koszul cohomology groups of surfaces and
higher dimensional varieties -- as Green put it,  there are more
reasons to compute them than  ways of computing them. In our
experience, this is especially so if the adjoint linear series
involves base-point-free or ample line bundles, to mention the
case of Mukai's conjecture. In this article and in previous ones we
have developed techniques to compute these cohomology groups.
Firstly in the proofs of the  vanishings leading to results on
higher syzygies we use induction on the number of ample and
base-point-free line bundles composing the line bundle we are
studying. 
To prove the vanishings which correspond to the first step of the
induction we have found it necessary to use the intrinsic 
geometric properties of the varieties under consideration. We make
here a distinction between two classes of varieties: those with
irregularity 
$h^1(\Cal O_X) > 0$ and those with irregularity $h^1(\Cal O_X)
= 0$. In the former case we use arguments involving 
Castelnuovo-Mumford regularity and the existence of
enough homologically trivial line bundles to show the
surjectivity of certain multiplication maps of vector bundles on
the variety.  In the latter case (comprising among others  K3
surfaces, Fano varieties,  Enriques surfaces, anticanonical
rational surfaces and Calabi-Yau
threefolds) we give  uniform proofs using induction on the
dimension of  the variety. Precisely we choose a suitable divisor
(we  point out to the reader that
it is not a hyperplane section!) on the variety to reduce the
question of the surjectivity of multiplication maps on the ambient
variety to a question of  surjectivity of multiplication maps
on the divisor.  This allows us to use eventually semistability
results and results on surjectivity of multiplication maps of
vector bundles on curves, like the technical (and beautiful)
results by Butler and Pareschi, 
\Butler \ns , Proposition 2.2 and \Pareschitwo \ns,
Corollary 4. 
These methods just introduced are 
displayed in more detail in the first sections, especially in Section
1 
 and Section 2. Later on  similar arguments are dealt  with sometimes
in a less detailed way.

\vskip .2 truecm

\noindent {\bf Convention.} Throughout this article we
work over an algebraically closed field   
of characteristic
$0$. For us  surface will  always mean minimal and smooth
algebraic surface. We will denote numerical equivalence of line 
bundles
by $\equiv$.

\vskip .2 truecm
\noindent{\bf Definition.} Let $X$ be a projective variety and let $L$
be a very ample line bundle on $X$. We say that $L$ is normally
generated or that satisfies the property $N_0$, if $|L|$ embeds $X$ as
a projectively normal variety. We say that $L$ is normally presented or
that $L$ satisfies the property $N_1$ if $L$ satisfies property $N_0$
and, in addition, the homogeneous ideal of the image of $X$ by $|L|$ is
generated by quadratic equations. We say that $L$ satisfies the
property $N_p$ for $p >1$, if $L$ satisfies property $N_1$ and the free
resolution of the homogeneous ideal of $X$ embedded by $|L|$ is linear
until the $p$th-stage.

\heading  1. A general result on syzygies of algebraic varieties
\endheading 
As we mentioned in the introduction, Green
interpreted the Betti numbers of the minimal free resolution of
the coordinate ring of an embedded projective variety in terms of
Koszul cohomology. Concretely, let
$X$ be a projective variety, and  let $F$ be a globally generated
vector bundle on $X$. We define the bundle $M_F$ as follows:
$$ 0 \to M_F \to H^0(F) \otimes \Cal O_X \to F \to 0 \ .\leqno
\sequence
$$
If $L$ is an ample line bundle on $X$ and all its positive powers are
nonspecial one has the following characterization of the property
$N_p$:

\proclaim {\GLlemma } Let $L$ be an ample, globally generated line
bundle on a variety $X$. If  
\text{$H ^1(\bigwedge ^{p'+1} M_L \otimes
L^{\otimes s})$} vanishes for all $0 \leq p' \leq p$ and all $s \geq
1$, then
$L$ satisfies the
 property
$N_p$. If in addition  $H^1(L ^{\otimes r}) = 0$, for all $r \geq 1$,
then the above  is a necessary and sufficient condition
for $L$ to satisfy property $N_p$. 
\endproclaim

\vskip .2 truecm
We will obtain our results on syzygies using the
previous lemma. For the proof of it we refer to
\EL
\ns, Section 1. Recall that we are working over an
algebraically closed field of characteristic $0$, thus in
our proofs we will check the vanishings of  
\text{$H ^1( M_L^{\otimes p'+1} \otimes
L^{\otimes s})$} rather than see directly the vanishings of 
\text{$H ^1(\bigwedge ^{p'+1} M_L \otimes
L^{\otimes s})$}. 

The purpose of this section is to prove a general result about 
Koszul cohomology and, by the above lemma, about  syzygies of
varieties of arbitrary dimension.

\proclaim{\preg}
Let $X$ be a projective variety. Let $B$ be a base-point-free line 
bundle on $X$ with regularity $r$. If $n\geq \text{ max}(r+p-2,p)$,
$p
\geq 1$ and $m\geq \text{ max}(r,1)$, then
$$H^i(M_{B^{\otimes m}}^{\otimes p+1} \otimes B^{\otimes
n+2-i})=0
\text{ for all }i 
\geq 1 \ .$$
In particular, $H^1(M_{B^{\otimes m}}^{\otimes p+1}\otimes
B^{\otimes n+1})=0$  
and if $B$ is ample and $n \geq \text{ max}(r+p-2,r,p)$,
then  $B^{\otimes n+1}$ satisfies the property $N_p$.
\endproclaim

To  prove the theorem we will need the following
\proclaim{\preglemma}
Let $X$ be and $B$ be as in  \preg\ns. If $n\geq r-1$ and $m\geq 1$,
then 
$$
H^{1}(M_{ B^{\otimes m} }\otimes B^{\otimes n+1})=0 \ . $$
In particular, if $B$ is ample,
then  $B^{\otimes n+1}$ satisfies the property $N_0$.
\endproclaim

\noindent {\it Proof.}
Since $n+1 \geq r$, $H^{1}( B^{\otimes n+1})=0$. Thus, 
tensoring the sequence \sequence relative to $B^{\otimes m}$ with $
B^{\otimes n+1}$ and taking global sections one sees that it is enough
to check that the multiplication map 
$$
H^{0}( B^{\otimes m})\otimes H^{0}( B^{\otimes n+1})\longrightarrow 
H^{0}( B^{\otimes m+n+1})
$$
surjects. To see that, we use the following useful observation:
\vskip .1 truecm
\noindent {\bf \sfreeobs \ns.}
Let $E$ and $L_1, \dots , L_r$ be 
coherent sheaves on a variety 
$X.$
Consider the map $H^0(E) \otimes H^0(L_1
\otimes \cdots \otimes L_r) @> \psi >> H^0(E \otimes
L_1 \otimes \cdots \otimes L_r)$ and the maps 
$$\displaylines{H^0(E) \otimes H^0(L_1
) @> \alpha_1 >> H^0(E \otimes L_1),
\cr
H^0(E\otimes L_1) \otimes H^0(L_2
) @> \alpha_2 >> H^0(E \otimes L_1\otimes L_2),\cr
 \dots ,
\cr
 H^0(E \otimes L_1 \otimes \cdots \otimes L_{r-1}) \otimes
H^0(L_r) @> \alpha_r >> H^0(E \otimes
L_1 \otimes \cdots \otimes L_r) \ .}$$
If $\alpha_1, \dots , \alpha_r$ are surjective
then $\psi$ is also surjective.
\vskip .1 truecm

In our case, we set $L_i=B$ and $E=B^{\otimes n+1}$,  and to see
that the maps
$\alpha_i$ are surjective we use the following 
generalization by
 Mumford
of a lemma of Castelnuovo (see \Mumford \ns; note that
the assumption of ampleness is unnecessary): 
\vskip .1 truecm
\noindent {\bf \CM \ns.}
 Let $L$ be a base-point-free line bundle
on a variety $X$ and let $\Cal F$ be a coherent sheaf on $X$. If
$\text H^i (\Cal F \otimes L ^{-i}) = 0$ for all $i \geq 1$, then
the multiplication map 
$$ \text H ^0(\Cal F \otimes L ^{\otimes i}) \ \otimes \text
H ^0(L) \to  \text H ^0(\Cal F \otimes L ^{\otimes i+1})$$
is surjective for all $i \geq 0$.
\vskip .1 truecm

Finally, the vanishings  required according to (1.4.2)
follow from our  assumption on regularity.
$\square$

\vskip .2 truecm

\noindent (1.5) {\it Proof of \preg \ns.}
The proof is by induction on $p$. We prove the result for $p=1$. 
First we show  that
$$
H^{1}(M_{B^{\otimes m}}^{\otimes 2}\otimes 
B^{\otimes n+1})=0\text{ for all }m\geq
r,1 \text{ and all } n \geq r-1,1
\ .$$
We will use \CM and \sfreeobs to prove this
statement. Observe that tensoring the sequence \sequence 
with $
M_{ B^{\otimes m}}\otimes B^{\otimes n+1}$ and taking global 
sections yields the following
long exact sequence:
$$\displaylines{
H^{0}(M_{ B^{\otimes m} }
\otimes B^{\otimes n+1})\otimes H^{0}( B^{\otimes m})
@>\gamma >>  
 H^{0}(M_{ B^{\otimes m}}\otimes B^{\otimes m+n+1})
\cr
\longrightarrow H^{1}(M_{ B^{\otimes m}}^{\otimes 2}
\otimes B^{\otimes n+1})
\longrightarrow
H^{1}(M_{ B^{\otimes m}}\otimes B^{\otimes n+1})
\otimes H^{0}( B^{\otimes m}).}$$
The last term in the above sequence is zero by 
\preglemma . 
Thus it is
enough to prove that $\gamma $ surjects. 
By \sfreeobs it is enough  to show that the multiplication
map 
$$
H^{0}(M_{ B^{\otimes m}}\otimes B^{\otimes n+1})\otimes H^{0}(B)
\longrightarrow
H^{0}(M_{ B^{\otimes m}}\otimes B^{\otimes n+2})
$$
surjects for all $m \geq r,1$ and all $n \geq r-1,1$. 
Since $B$ is base-point-free, by \CM we need to check the
 vanishings
$ H^i(M_{ B^{\otimes m}}\otimes B^{\otimes n+1-i})=0
\text{ for all  }  i \geq 1$, l $m \geq r,1$ and $n \geq r-1,1$.
For $i \geq 2$, we tensor the sequence \sequence corresponding to 
$ B^{\otimes m}$ with $ B^{\otimes n+1-i}$ and take global
sections.  The vanishings then follow from our assumption on the 
regularity of $B$. Since $m \geq r$ and $n \geq r-1$ it
follows in particular that
$H^1(B^{\otimes m})=H^1(B^{\otimes n})=0$, hence the
vanishing required  for $i=1$ is
equivalent to the vanishing of $H^1(M_{B^{\otimes n} } 
\otimes
B^{\otimes m})$, which follows in turn from \preglemma 
\ns.

The vanishings of $
H^{i}(M_{B^{\otimes m}}^{\otimes 2}\otimes B^{\otimes n+2-i}) 
\text{ for all }m\geq
1, \text{ all } i \geq 2 \text{ and all } n \geq r-1
$
follow from \sequence \ns, \preglemma \ns, and the assumption on 
regularity.

Let us now assume that the desired vanishings occur for $p-1$.  
We therefore have:
$$
\displaylines{H^{i}(M_{B^{\otimes m}}^{\otimes p}\otimes
B^{\otimes n+2-i})=0 
\text{ for all }n
\geq \text{ max}(p+r-3,p-1), \cr 
\text{ all } m \geq \text{
max}(r,1) \text{ and all } i
\geq 1
\ .}
$$
We first prove the desired vanishing for $p$ and $i=1$. 
By tensoring 
the sequence \sequence with 
$M_{B^{\otimes m}}^{\otimes p}
\otimes B^{\otimes n+1}$
and taking global sections one sees that the desired 
vanishing can be
 obtained by
showing the surjectivity of the multiplication map $\delta $ sitting in
the long exact sequence 
$$\displaylines{
H^{0}( M_{B^{\otimes m}}^{\otimes p}\otimes B^{\otimes n+1})\otimes
H^{0}(B^{\otimes m})
@> \delta >> H^{0}( M_{B^{\otimes m}}^{\otimes p}\otimes
B^{\otimes m+n+1})
\cr
\longrightarrow H^{1}( M_{B^{\otimes m}}^{\otimes p+1}\otimes 
B^{\otimes
n+1})\longrightarrow H^{1}( M_{B^{\otimes m}}^{\otimes p}
\otimes B^{\otimes n+1})
\otimes H^{0}(B^{\otimes m}) \ .} 
$$
The last term is zero by induction assumption. In order to
prove the surjectivity of $\delta $ we use \sfreeobs \ns. 
By \sfreeobs it suffices to  show
the
surjectivity of the  map

\noindent $
H^0(M_{B^{\otimes m}}^{\otimes p}
\otimes B^{\otimes n+1})
\otimes H^0(B) @> 
{\epsilon
} >> H^0(M_{B^{\otimes m}}^
{\otimes p}\otimes B^{\otimes n+2})$  for all $n \geq p+r-2, p$
and all $m \geq r, 1$ .

\noindent To prove the surjectivity of $\epsilon$ we use \CM \ns. 
According to  it, it suffices 
that the groups 
$H^i(M_{B^{\otimes m}}^{\otimes p}\otimes B^{\otimes
n+1-i})$ vanish, which  follows by
induction.

Finally, to show that
$H^i(M_{B^{\otimes m}}^{\otimes p+1}\otimes B^{\otimes
n+2-i})=0 
\ , \text{ for all } i \geq 2$ we consider again sequence
\sequence associated to $B^{\otimes m}$,  tensor it with
$M_{B^{\otimes m}} ^{\otimes p}\otimes B^{\otimes n+2-i}$
and take global sections. Then the vanishings follow again
from induction  hypothesis.

The fact that $B^{\otimes n+1}$ satisfies the property
$N_p$ follows from the vanishing of \linebreak
$H^1(M_{B^{\otimes
n+1}}^{\otimes p'}
\otimes B^{\otimes s(n+1)})$ for all $1 \leq p' \leq p$ and all $s \geq 1$,
from \preglemma and from \GLlemma \ns.
$\square$
\vskip .2 truecm

The theorem just proven, which might seem at first glance somehow
vague,  holds however the power to unify several results for
different kinds of varieties: it yields information about 
pluricanonical
embeddings of surfaces of general type (Corollary 5.11). It also
implies the following corollary  concerning varieties of arbitrary
dimension and canonical divisor numerically trivial, an
 example of which are Calabi-Yau n-folds:

\proclaim{Corollary 1.6}
Let $X$ be a variety of dimension $m$ with $K_X \equiv 0$ and let 
$B$ be ample and base-point-free line bundle. Let $L= B^{\otimes
n+1}$. If $n \geq p+m-1$, then $L$ satisfies property $N_p$. In
particular, if $X$ is a Calabi-Yau threefold, $B$ is
an ample and base-point-free line bundle on $X$, $n \geq p+3$
and $p
\geq 1$,  then
$B^{\otimes n}$ satisfies property $N_p$.
\endproclaim

\noindent {\it Proof.} The result is a straight forward
consequence of 
\preg \ns, since by Kodaira vanishing Theorem, $B$ is
$(n+1)$-regular. $\square$
\vskip .2 truecm

\preg also implies a result for surfaces with $p_g =0$
(among them  elliptic ruled surfaces, Enriques surfaces
and bielliptic surfaces):
\proclaim{\GPtwo \ns, Theorem 2.2} Let
$X$ be  a surface with
$p_g = 0$. 
Let
$B$ be a nonspecial, ample, and base-point-free line
 bundle. Then
$B^{\otimes p +1}$ satisfies the property $N_p$  for all 
$p \geq
1$.
\endproclaim
Therefore \preg and its corollaries are a good starting point for our
study of syzygies of varieties. However, if one focuses on the
particular examples and uses the specific geometry of the
varieties in question, one can expect to obtain sharper and more
complete results. Precisely this was done for elliptic ruled surfaces
in
\GPtwo and is done for Enriques surfaces in Section 2, for
bielliptic surfaces in Section 4, for surfaces of general type
in Section 5.

\heading 2. Syzygies of Enriques surfaces \endheading

In Section 1 we proved a general theorem, \preg \ns, which unifies
a number of results for different kinds of varieties. In this
section we focus on Enriques surfaces. The geometric genus of an
Enriques surface is $0$ and, in characteristic
$0$, a globally generated line bundle over an Enriques surface has
null higher cohomology, hence it is $2$-regular. Therefore the
starting point of our study of syzygies of Enriques surfaces is the
following theorem, corollary of \preg \ns, which 
fits indeed in \text{\MP \ns:}

\proclaim{\Enriquessame \ns, 
(cf. \GPtwo \ns, Corollary 2.7.1)} 
Let $X$ be an Enriques surface. Let
$B$ be a base-point-free line bundle. Then the
image of $X$ by $|B^{\otimes p+1}|$ satisfies
property $N_p$, for all
$p \geq 1$. If in addition $B$ is ample then
$B^{\otimes p+1}$ is very ample and satisfies the property
$N_p$, for all
$p \geq 1$.
\endproclaim

Our intention now is to study a more general class of line bundles
(namely, tensor products of $p+1$ different base-point-free line
bundles),
and in particular, adjoint line bundles. For that we need to follow a
different approach: roughly, we are going to use ``induction on the
dimension", in the sense explained in the introduction. This approach
will unfold throughout this section and the machinery developed along
the way will be used for other results of this article, concretely in
Sections 3 and 5. We now resume with a result about normal generation:

\proclaim{\Enriquesng} Let $X$ be 
an Enriques surface. Let $B_1$, $B_1'$, $B_2$ and
$B_2'$ be  
ample and base-point-free line bundles on
$X$,
such that $B_1 \equiv B_1'$, $B_2 \equiv B_2'$ and,
either
$B_1 \cdot B_2 \geq 4$,  
$B_1^2
\geq 6$, and $B_2^2
\geq 6$ or $B_1 \cdot B_2 \geq 5$. Let
$L=B_1^{\otimes r}
\otimes B_2^{\otimes s}$ and $L'={B_1'}^{\otimes
k}\otimes {B_2'}^{\otimes l}$. If $r,s ,k \geq 1, \text {
and } l\geq 0$, then the map  
$H^0(L)\otimes H^0(L') @> \alpha >> H^0(L\otimes L')$ surjects and
$H^1(M_L \otimes L')=H^1(M_L'\otimes L)=0$. In particular, $L$ is very
ample and satisfies property $N_0$. 
\endproclaim 

Before we go on with the proof of \Enriquesng \ns , we isolate for
convenience three ingredients of the argument, which will
be used in many other instances. 
The first is an observation on the relation between 
the surjectivity of
multiplication maps, and the surjectivity of its restrictions
to divisors. The other two  are a result due to Butler
and another one due to Pareschi,  about the surjectivity
of multiplication maps of vector bundles on curves.

\proclaim{\robs} Let $X$ be a regular variety ( i.e, a variety such that 
$H^1(\Cal{O}
_X)=0).$ Let $E$ be a vector bundle on $X$, let $C$ be a
divisor such that  
$L$
$=\Cal{O}_X\left( C\right) $ is  globally generated line bundle and 
$H^1(E\otimes L^{-1})=0.$ If the multiplication map 
$
H^0(E\otimes \Cal{O}_C)\otimes H^0(L\otimes
\Cal{O}_C)\to H^0(E\otimes L\otimes \Cal{O}_C)
$
surjects, then the map 
$
H^0(E)\otimes H^0(L) \to H^0(E\otimes L)
$
also surjects.
\endproclaim

\noindent {\it Proof.} We construct the following
commutative diagram:
$$\matrix
H^0(E) \otimes H^0(\Cal O_X)&\hookrightarrow &H^0(E)
\otimes H^0(L)&\twoheadrightarrow &H^0(E) \otimes H^0(L
\otimes 
\Cal O_C)\cr
\downarrow&&
\downarrow&&\downarrow\cr
H^0(E )&\hookrightarrow &H^0(E \otimes
L)&\twoheadrightarrow &H^0(E \otimes L \otimes \Cal
O_C)\ .
\endmatrix
$$
The surjectivity of the left hand side vertical map is
obvious. The surjectivity of the right hand side
vertical map follows by hypothesis. The exactness of the
top horizontal sequence follows from the fact that $X$
is regular. The claim is the surjectivity of the middle
vertical map. $\square$

\proclaim{\Butres  (\Butler
\ns , Proposition 2.2)} Let
$E$ and
$F$ be semistable vector bundles over a curve $C$ such that $E$ is
generated by its global sections. If
\roster
\item $\mu(F)>2g$, and
\item $\mu(F)>2g+\text{ rank}(E)(2g-\mu(E))-2h^1(E)$,
\endroster
then the multiplication map $H^0(E) \otimes H^0(F) \to 
H^0(E\otimes
F)$ surjects.
\endproclaim

\proclaim{\Pareschilemma \ns (\Pareschitwo \ns \ns ,
Corollary 4.)}
Let $N$ and $L$ be two base-point-free line bundles on
$C$ such that:
\roster
\item"(a)" at least one of them is very ample;
\item"(b)" $h^0(N), h^0(L) \geq 3$ and 
\item"(c)" deg$N+ \text{ deg }L \geq \text{ max }(3g-3,
4g+1-2h^1(N)-2h^1(L) -\text{Cliff}(C))$.
\endroster
Then the multiplication map
$$H^0(L) \otimes H^0(N) \longrightarrow H^0(L \otimes
N)$$
is surjective.
\endproclaim

\noindent (2.6) {\it Proof of \Enriquesng \ns .}
Note first that, since we are working over a field of
characteristic $0$, any base-point-free line bundle on
$X$ has null higher cohomology. If we twist the
sequences \sequence relative to $L$ and $L'$ by $L'$ and
$L$ respectively and take global sections, we see at once that
$H^1(M_L \otimes L') = H^1(M_L' \otimes L)$ and  equal to the
cokernel of
$$H^0(L) \otimes H^0(L') @> \alpha >> H^0(L \otimes L') \ .$$ To see that
$\alpha$ indeed surjects, we use \sfreeobs \ns . According 
to it we
want to check that several (possibly more than one) multiplication
maps surject. We check here the first one; the surjectivity of the rest
can be seen in the same way. The  map in question is
$$H^0(L) \otimes H^0(B_1') @> \beta >> H^0(L \otimes
B_1') \ .$$ To see the surjectivity of $\beta$, we
consider a smooth irreducible curve
$C$  in
$|B_1'|$ (such curve exists by Bertini's Theorem because
$B_1'$ is ample and base-point-free) and   use
\robs \ns. It is therefore enough to check that
$$H^0(L \otimes \Cal O_C) \otimes H^0(B_1' \otimes \Cal
O_C) @>
\gamma >> H^0(L
\otimes B_1' \otimes \Cal O_C) $$
surjects. For that, if $B_1 \cdot B_2 \geq 5$, we may apply \Butres
\ns . Indeed, the line bundle $B_1'$ is globally
generated, and by adjunction $\mu (L) = \text{ deg}L
\geq 2g(C)+3 > 2g(C) + 2$. If $B_1
\cdot B_2 = 4$ and $B_1^2 \geq 6$, then $g(C) \geq 4$ and, since
$C$ is irreducible, it follows that it is non-hyperelliptic (cf. \CosDol  \ns,
Proposition 4.5.1). Then  the surjectivity of $\gamma$ follows
from \Pareschilemma \ns.
$\square$

\vskip .2 truecm
As a corollary of \Enriquesng we prove a stronger
version of the conjecture of Mukai, in the case of
Enriques surfaces and for the property
$N_0$. To see that we use the following

\proclaim {\twoamplearefree}
Let $A_1$  and $A_2$ be two ample divisors on a surface $X$ with
Kodaira dimension 0. Then 
$ A_1 \otimes A_2$ is base-point-free.
\endproclaim

\noindent {\it Proof.}
Since $K_X \equiv 0$,  $(A_1 \otimes A_2)^2\geq 5$.  By hypothesis
$A_1 \otimes A_2$ is ample. If $A_1 \otimes A_2$
were not base-point-free, it would follow from Reider's theorem that
there would exist an
effective divisor $E$ such that one of the following holds:
\roster

\item"(a)"$(A_1 \otimes A_2) \cdot
E=0$
and $E^{2}=-1$ or 

\item"(b)" $(A_1 \otimes A_2) \cdot E=1$ and $E^{2}=0.$ 
\endroster
None of the two possibilities can occur since, $A_i$ being ample, $A_i
\cdot E \geq 1$.
$\square$
\vskip .2 truecm

\proclaim{\EnriquesMukaing}
Let $X$ be an Enriques surface and $A_1, \dots , A_n$   ample line
bundles on
$X$. Let $L = K_X \otimes A_1 \otimes \cdots \otimes A_n$. If $n \geq
4$, then
$L$ satisfies property $N_0$.
\endproclaim

\noindent {\it Proof.}
By \twoamplearefree \ns, $K_X \otimes A_1 \otimes A_2$ and $A_3
\otimes
\cdots \otimes A_n$ are base-point-free line bundles.
There are furthermore ample, and, by adjunction, $(K
\otimes A_1 \otimes A_2)^2 \geq 6$, $(A_3 \otimes \cdots 
\otimes A_n)^2
\geq 6$ and 
$(K
\otimes A_1 \otimes A_2)\cdot (A_3 \otimes \cdots
\otimes A_n)
\geq 4$. Then the result follows from \Enriquesng \ns.
$\square$
\vskip .2 truecm

We now generalize these results to
higher syzygies. To do so, we need another two lemmas. 
In the case in which $\frak q$ is a curve $C$, the former
allows us to pass from a multiplication map involving
non-semistable bundles (note 
that $M_F \otimes \Cal O_C$ is unstable if
$H^1(L \otimes \Cal O(-C))=0)$ to a multiplication map
involving  semistable bundles.  This  situation  is of
course easier to handle. 
The latter lemma deals with positivity and semistability of
bundles on curves.
They will not only be used for the arguments in the
remaining of this section but also in Section 3 and 5.

\proclaim{\Splemma } Let $X$ be a projective variety, let $q$ be a
nonnegative integer and let $F_i$ be a base-point-free line bundle
on $X$ for all $1 \leq i \leq q$. Let $Q$ be an effective line
bundle on $X$ and let $\frak q$ be a reduced and irreducible member
of
$|Q|$. 
Let $R$ be a line bundle and $G$ a sheaf on $X$ such that
\vskip .1 cm
\item{1.} $\text H^1(F_i \otimes Q^*)=0$
  \item {2.}  
$
\text H
^0(M_{(F_{i_1}
\otimes
\Cal O_\frak q)}
\otimes \dots \otimes M_{(F_{i_{q'}} \otimes \Cal O_\frak q) }\otimes
R \otimes \Cal O_\frak q)
\otimes \text H^0(G) \to$
\item{}
$\to \text H ^0(M_{(F_{i_1} \otimes 
\Cal O_\frak
q)}
\otimes \dots \otimes M_{(F_{i_{q'}}  \otimes \Cal O_\frak
q)}
\otimes R \otimes G \otimes \Cal O_\frak q) \ \text {surjects for all}
\  0 \leq q' \leq q$.

\noindent Then, for all $0 \leq q'' \leq q$ and any subset
$\{j_k\}
\subseteq \{i\}$ with $\#\{j_k\}=q''$ and for all $0 \leq k'
\leq q''$, 
$$\displaylines {\text H^0(M_{F_{j_1}} \otimes \dots \otimes
M_{F_{j_{k'}}}
\otimes M_{(F_{j_{k'+1}} \otimes \Cal O_\frak q)} \otimes \dots
\otimes M_{(F_{j_{q''}}\otimes \Cal O_\frak q)} \otimes R
\otimes \Cal O_\frak q) \otimes \text H^0(G) \to \cr
\text
H^0(M_{F_{j_1}} 
\otimes \dots \otimes M_{F_{j_{k'}}}
\otimes M_{(F_{j_{k'+1}} \otimes \Cal O_\frak q)} \otimes \dots
\otimes M_{(F_{j_{q''}} \otimes \Cal O_\frak q)} \otimes G \otimes R
\otimes \Cal O_\frak q)}$$ surjects.
\endproclaim

\noindent {\it Proof.}
We prove the result by induction on $q''$. For $q''=0$ the
corresponding statement is just Condition 2 when $q=0$. 
Assume that
the result is true for $q''-1$. In order to prove the result
for $q''$ we will use induction on $k'$. If $k'=0$, the
statement is again just Condition 2. Assume that the result is true
for
$k' -1$. Now for any $F$ globally generated vector bundle and
for any effective divisor $\frak q$ such that $H^1(F \otimes Q^*)= 0$,
for  $Q=
\Cal O(\frak q)$, we have this commutative diagram: 
$$\matrix
&&0&&0&& \\
&& \downarrow && \downarrow && \\
0 &\to& \text H ^0(F \otimes Q ^*) \otimes \Cal O _\frak q &
\to & \text H ^0(F \otimes Q ^*) \otimes \Cal O _\frak q 
&\to& 0 \\
&& \downarrow && \downarrow && \downarrow \\
0& \to & M_{F} \otimes \Cal O _\frak q & \to  & \text H ^0(F
) \otimes \Cal O_\frak q 
&\to&F \otimes \Cal O_\frak q & \to & 0 \\
&& \downarrow && \downarrow && \downarrow \\
0& \to& M_{(F \otimes \Cal O _\frak q)}&\to&\text H ^0(F
\otimes
\Cal O _\frak q) \otimes \Cal O _\frak q& \to &F \otimes
\Cal O _\frak q &\to &0 \\
&& \downarrow && \downarrow && \downarrow \\
&&0&&0&&0 \\
\endmatrix$$
We are interested in the left hand side vertical exact
sequence:
$$  0 \to  \text H ^0(F \otimes Q ^*) \otimes \Cal O
_\frak q \to M_{F} \otimes \Cal O _\frak q \to 
M_{(F \otimes \Cal O _\frak q)} \to 0 \leqno \Splemmaone $$
By
Condition 1,
$F$ can be taken to be
$F_{j_{k'}}$. Tensoring \Splemmaone by
$$M_{F_{j_1}} \otimes \dots
\otimes M_{F_{j_{k'-1}}}
\otimes M_{(F_{j_{k'+1}} \otimes \Cal O_\frak q)} \otimes \dots
\otimes M_{(F_{j_{q''}}\otimes \Cal O_\frak q)} \otimes R
\otimes \Cal O_\frak q \ ,$$ taking global sections and tensoring by
$\text H^0(G)$ we obtain this  commutative diagram:
$$ \matrix 0 & \to & A \otimes \text H^0(G)& \to & B \otimes \text
H^0(G) &
\to & C
\otimes \text H^0(G) & \to & 0 \\
&&\downarrow & & \downarrow & & \downarrow \\
0 & \to & A' & \to & B' & \to & C' & & \\
\endmatrix
$$
where $A = \text H^0(F_{j_{k'}} \otimes Q^*) \otimes 
\text H^0(\bigotimes_{r=1}^{k'-1}M_{F_{j_r}} \otimes 
\bigotimes_{r=k'+1}^{q''} M_{(F_{j_r} \otimes \Cal O_\frak q)} 
\otimes R
\otimes \Cal O_\frak q)$, 
 $B = \text 
H^0(\bigotimes_{r=1}^{k'}M_{F_{j_r}} \otimes 
\bigotimes_{r=k'+1}^{q''} M_{(F_{j_r} \otimes \Cal O_\frak q)} 
\otimes R
\otimes \Cal O_\frak q)$,  
$C = \text
H^0(\bigotimes_{r=1}^{k'-1}M_{F_{j_r}} \otimes 
\bigotimes_{r=k'}^{q''} M_{(F_{j_r} \otimes \Cal O_\frak q)} 
\otimes R
\otimes \Cal O_\frak q)$, 
$A' = \text H^0(F_{j_{k'}} \otimes
Q^*)
\otimes 
\text H^0(\bigotimes_{r=1}^{k'-1}M_{F_{j_r}} \otimes 
\bigotimes_{r=k'+1}^{q''} M_{(F_{j_r} \otimes \Cal O_\frak q)} 
\otimes R
\otimes G \otimes \Cal O_\frak q)$, $B' = \text 
H^0(\bigotimes_{r=1}^{k'}M_{F_{j_r}} \otimes 
\bigotimes_{r=k'+1}^{q''} M_{(F_{j_r} \otimes \Cal O_\frak q)} 
\otimes R
\otimes G \otimes \Cal O_\frak q)$ and $C' = \text
H^0(\bigotimes_{r=1}^{k'-1}M_{F_{j_r}} \otimes 
\bigotimes_{r=k'}^{q''} M_{(F_{j_r} \otimes \Cal O_\frak q)} 
\otimes R
\otimes G \otimes \Cal O_\frak q)$. The top horizontal exact sequence
is certainly surjective: this follows from chasing the diagram after
having taken cohomology. The left hand side vertical sequence surjects 
by the induction hypothesis on
$q''$ and  the right hand side exact sequence surjects by induction
on $k'$ (we have assumed the result to be true for $q''-1$ and
$k'-1$). Therefore we obtain the surjectivity of the vertical
sequence sitting  in
the middle  of
the commutative diagram. \ $\square$

\proclaim{\slopelemma}
Let $E$  be a semistable 
vector bundle with  
$\mu (E) > 2g $ and $F$ a vector bundle on a curve $C$ of
genus $g$. 
\roster
\item"(1)" If 
 $\mu (F) \geq 2g+4 $, then $\mu
(M_{E}\otimes F)> 2g+2$. 
\item"(2)" If $\mu (F) \geq 2g+2 $, then $\mu
(M_{E}\otimes F)> 2g$.
\endroster
Moreover, if $F$ is in addition
semistable, then 
$M_{E}\otimes F$ is semistable.
\endproclaim

\noindent {\it Proof.}
 Since $E$ is semistable and $\mu(E) >2g $, $E$ is globally
generated and $h^1(E)=0$, hence the vector bundle $M_E$ is
defined and has slope 
$$\mu (M_E)=\frac {-\mu(E)} {\mu(E)-g} \ .$$
Then, for (1),  $\mu(M_E\otimes F)
\geq \frac {-\mu(E)} {\mu(E)-g} + 2g+4$. Thus if $\frac {-\mu(E)}
{\mu(E)-g} + 2g+4 > 2g+2$ we are done, but that inequality is
equivalent to $\mu(E) > 2g$. The proof of (2) is analogous. 
Now, if $F$ is semistable by \Butler \ns, Theorem 1.12 and 
\Miyaoka
\ns, Corollary 3.7, 
$M_E
\otimes F$ is also semistable.
$\square $

\proclaim{\Enriquesnpres}
Let $X$ be an Enriques surface. Let $B_{1}$, $B_1'$,
$B_2$ and B$_{2}'$  be
two  ample and 
base-point-free divisors such that $B_1 \equiv B_1'$,
$B_2 \equiv B_2'$ and 
$B_{1}\cdot B_{2}\geq 6.$ Let $L= B_1^{\otimes s} 
\otimes B_2^{\otimes r}$
and $L' = {B_1'}^{\otimes k} \otimes {B_2'}^{\otimes
l}$.  If $k,l,r,s \geq 1$,
then 
$
H^{1}(M_L^{\otimes 2}\otimes
L')=0
$. In particular, 
$L'$ satisfies property $N_1$.
\endproclaim

\noindent{\it Proof.}
 The cohomology group 
$H^{1}(M_L^{\otimes 2}\otimes
L')=0$
sits in the  long exact
sequence 
$$\displaylines{
H^{0}(L)\otimes H^{0}(M_{L}\otimes L^{\prime })@>\alpha
>> H^{0}(M_{L}\otimes L\otimes L^{\prime }) \cr
\longrightarrow H^{1}(M_{L}^{\otimes 2}\otimes L^{\prime
})\longrightarrow
H^{0}(L)\otimes H^{1}(M_{L}\otimes L^{\prime })
\ ,}$$
obtained by tensoring \sequence relative to $L$ with
$M_L\otimes 
L'$    and taking global sections.
The last term is zero by \Enriquesng \ns,  
thus it is enough to prove that
$\alpha $ is
surjective. To show the surjectivity of $\alpha $ we use 
\sfreeobs \ns.
According to it we need to check the  surjectivity of several
maps. Here we will only show the
surjectivity of the first of them, since the rest are
analogous:
$$
H^{0}(B_{1})\otimes H^{0}(M_{L}\otimes L^{\prime })@>\beta >>
H^{0}(M_{L}\otimes L^{\prime }\otimes B_{1})
\ .$$
Let $C$ be a smooth member of $|B_1|$.
From \Enriquesng it follows that  
$H^{1}(M_L\otimes L' \otimes B_1^*)=0$, 
therefore we may apply
\robs to reduce the question of surjectivity of 
$\beta $ to the
surjectivity of the following multiplication map on 
$C$: 
$$
 H^{0}(B_{1}\otimes
\Cal{O}_{C})\otimes H^{0}(M_{L}\otimes L^{\prime
}\otimes \Cal{O}_{C}) \longrightarrow
H^{0}(M_{L}\otimes L^{\prime }\otimes B_{1}\otimes \Cal{O}_{C})
\ . $$
 By
\Splemma it is enough to check that the following 
multiplication maps on $C$
are surjective: 
$$\displaylines{
H^0(B_1 \otimes \Cal O_C) \otimes H^0(L' \otimes \Cal
O_C) \longrightarrow H^0(B_1 \otimes L' \otimes \Cal
O_C) \cr
H^{0}(B_{1}\otimes \Cal{O}_{C})\otimes H^{0}(M_{L\otimes \Cal{O}
_{C}}\otimes L^{\prime }\otimes \Cal{O}_{C})@>\gamma  >>
H^{0}(M_{L\otimes \Cal{O}_{C}}\otimes L^{\prime
}\otimes B_{1}\otimes \Cal{O}_{C})
\ .}$$
The surjectivity of the first map was already seen
within the course of proving \Enriquesng \ns. For
$\gamma$, we use \Butres \ns. Since deg$(L \otimes \Cal
O_C)$ and  deg$(L' \otimes \Cal O_C)$ are both
greater than or equal to $2g+4$, it follows from \slopelemma that the
bundle
$M_{L\otimes
\Cal{O} _{C}}\otimes L^{\prime }\otimes \Cal{O}_{C}$ 
is semistable with
slope strictly bigger than $2g+2$.
Then it follows from \Butres
that $\gamma $ is surjective and we are done.
Now, since $L'$ is ample, it follows from the vanishing of
$H^1(M_{L'}^{\otimes 2}\otimes {L'}^{\otimes s})$ for
all $s \geq 1$,
\Enriquesng and 
\GLlemma \ns,
 that ${L'}$ satisfies property $N_1$.
$\square$
\vskip .2 truecm
We obtain the following corollary, which proves  Mukai's conjecture
(and when considering powers of the same ample bundle,
improves his bound), regarding property
$N_1$  for Enriques surfaces. 
\vskip .3 truecm

\proclaim{\EnriquesMukainpres} Let $X$ be an Enriques surface. Let
$A$,
$A_1, \dots, A_n$ be ample line bundles. Then the line bundles
$K_{X}\otimes A^{\otimes m}$ and $K_{X}\otimes A_1 \otimes \cdots
\otimes A_n$ satisfy
property
$N_{1}$ if $m \geq 4$ and $n \geq 5$ respectively.
\endproclaim

\noindent {\it Proof.}
For the former case, let $B_1=K_X \otimes A^{\otimes 2}$ and $B_2 =
A^{\otimes m-2}$. For the latter, let $B_1=K_X \otimes A_1
\otimes A_2
\otimes A_3$ and $B_2 = A_4 \otimes \cdots \otimes
A_n$. In both cases, $B_1$ and
$B_2$ are ample, and base-point-free by \twoamplearefree \ns.
Furthermore,
$B_1 \cdot B_2 \geq 6$, consequently the result follows
from
\Enriquesnpres \ns.
$\square$

\vskip .2 truecm
To finish this section we show a result for higher syzygies 
of adjoint
bundles.
Before that we state a useful lemma dealing with the numerical
nature of the property of base-point-freeness.

\proclaim {\numfree}
Let $X$ be a surface with nonnegative Kodaira dimension and let
$B$ be an ample and base-point-free line bundle such that $B^2 \geq
5$. If $B' \equiv B$, then $K_X \otimes B'$ is ample and
base-point-free. In particular, if $\kappa(X)=0$, $B'$ is ample and
base-point-free for all $B' \equiv B$.
\endproclaim

\noindent {\it Proof.}
The line bundle $B'$ is ample because ampleness 
is a numerical
condition and has self-intersection greater than or 
equal to $5$. If
$K_X \otimes B'$ has base points, by Reider's theorem there is an 
effective divisor
$E$ such that:
\roster
\item"(a)"
$B'\cdot E=0$ and $E^{2}=-1$
or 
\item"(b)" $B'\cdot E=1$
and $E^{2}=0.$ 
\endroster
The former cannot happen because $B'$ is ample. 
We will
also
rule out (b).  The divisor $E$ must be  irreducible 
and reduced because $B'$ is
ample and $B'\cdot E = 1$. On the other hand, the arithmetic genus of
$E$ is greater than  or equal to
$1$. Now $B\cdot E = B' \cdot E =1$ so
$h^{0}(B\otimes
\Cal{O}_{E}) \leq1$.  Since
$B$ is base-point-free, $E$ should be a smooth rational curve and this
is a contradiction. 
$\square$

\proclaim{\Enriqueshighersyzygies} Let $X$ be an 
Enriques surface.
Let
$B$ be an ample and base-point-free line bundle such
that 
$B^2
\geq 6$ and let $N, N'$ be line bundles numerically
equivalent to $0$  
(i.e., they are either trivial or equal to
$K_X$). 
Let
$L = B^{\otimes p+1+l}
\otimes N$,    
$L^{\prime }= B^{\otimes p+1+k} \otimes N'$ 
for
$p
\geq 1$. Then 
$H^{1}(M_{L}^{\otimes p+1}\otimes L^{\prime })$
  vanishes for all  $k,l
\geq 0$. In particular
$L$ satisfies property $N_p$.
\endproclaim

\noindent {\it Proof.} Since
$B^2
\geq 6$,  by
\numfree the line bundle 
$B \otimes N$ is also ample and base-point-free.  The proof
is by induction. The result is true for
$p=1$ by  \Enriquesnpres \ns. We 
assume now the result to be true for $p-1$.   
In particular we have $H^{1}(M_{L}^{\otimes p}\otimes L')
=0$. Tensoring the
sequence \sequence with $M_{L}^{\otimes p}\otimes L'$ 
and
taking global sections yields therefore the following 
long exact
sequence
$$
H^{0}(L)\otimes H^{0}(M_{L}^{\otimes p}
\otimes
L')@>\alpha >> H^{0}(M_{L}^{\otimes p}\otimes L \otimes
L')
\to H^{1}(M_{L}^{\otimes p+1}\otimes
L')\to 0 \ ,
$$
thus  it is enough to
prove that
the multiplication map $\alpha $ is surjective. 
Then by 
\sfreeobs it is enough to see the surjectivity of  
$$
H^{0}(B')\otimes H^{0}(M_{L}^{\otimes p}\otimes
B^{\otimes p+1+k} \otimes N')@>{\beta }>>
H^{0}(M_{L}^{\otimes p}\otimes  B^{\otimes p+1+k} \otimes
B'
\otimes N')
\ ,
$$
where $B'$ is either $B$ or $B \otimes N$. Now to complete the proof
one can argue in two ways. One of them is using \CM \ns. 
The path to follow is shown in the proof of \preg but we outline
here the steps to be taken. The first cohomology 
vanishing
required, 
$$   H^{1}(M_{L}^{\otimes
p}\otimes B^{\otimes p+1+k} \otimes N' \otimes {B'}^*)
\leqno{(2.14.1)}  $$  follows directly by induction. For
the second cohomology vanishing one may observe that,
after iteratively chasing the cohomology sequence, it
follows by induction,  from
\Enriquesng and Kodaira vanishing Theorem. The other way 
to argue is as for the surjectivity of
$\beta$ in the proof of
\Enriquesnpres \ns: one uses  \Splemma to
reduce the problem to checking the surjectivity of
multiplication maps on a curve. 

Finally  since $L$ is ample, \GLlemma
implies that $L$ satisfies $N_p$.
$\square$

\proclaim{\EnriquesMukaiNp} Let $X$ be an 
Enriques surface, let $A$ be 
an ample line bundle and $B$ an ample and
base-point-free line bundle on
$X$. If
$m
\geq p+1$, then
$K_X
\otimes B^{\otimes m}$ satisfies property $N_p$. If $n \geq
2p+2$, then $K_X
\otimes A^{\otimes n}$ satisfies property $N_p$.
\endproclaim

\noindent {\it Proof.}
The first statement is a straight forward consequence of
the theorem.  By
\twoamplearefree \ns, the line bundle
$A^{\otimes 2}$ is base-point-free, so if $n$ is even the second
statement  follows from the first. If $n$ is odd the result follows from
a slight variation of the argument in the proof of
\Enriqueshighersyzygies \ns: we break up $K_X
\otimes A^{\otimes n}$ as tensor product of $n-1$ copies of
$B=A^{\otimes 2}$ and $B'=A^{\otimes 3}$, which is
base-point-free by \twoamplearefree \ns. When applying
\sfreeobs we take the last map among the $\alpha_i$
to be precisely the map involving $B'$. The reader
can easily verify that the vanishings needed in order to apply
\CM follow by induction or, eventually, by Kodaira vanishing
Theorem. 
$\square$

\heading  3. Koszul rings of Enriques surfaces \endheading

We have devoted Section 2 to the study of syzygies of embeddings of
Enriques surfaces. We show in particular a result, \Enriquesnpres
\ns, about normal presentation of line bundles which were the tensor
product of two base-point-free line bundles. Recall that the normal
presentation property means that the homogeneous ideal of the
(projectively normal) variety is generated by forms of degree $2$.
As already pointed out in the introduction, an interesting
algebraic property that many normally presented rings have is
the Koszul property. 
There exist many significant examples: canonical rings of curves
(cf. \FV \ns, 
\PP \ns),
rings of
curves of degree greater than or equal to $2g+2$ (cf. \Butler \ns,
\GPone \ns), elliptic ruled surfaces (cf. \GPone \ns, Theorem 5.8) and
those line bundles on Enriques surfaces which are normally
presented according to \Enriquessame (cf. \GPone \ns, Corollary 5.7).
This section provides yet one more case in favor of this philosophy:
we will show in \EKoszul that those line bundles on an Enriques
surface which are normally presented according to 
\Enriquesnpres 
also
satisfy the Koszul property. Moreover, in the course of proving the
result, it can be seen how the property $N_1$ is one of the first
conditions required for the ring to be Koszul.

To begin we recall  some  notation and some basic
definitions:
given a line bundle
$L$ on a variety
$X$, we set $R(L) = \bigoplus _{n=0} ^\infty \text H
^0(X,L ^{\otimes n})$.

\proclaim {\defK} Let $R= \bold k \oplus R_1 \oplus R_2
\oplus
\dots$ be a graded ring and $\bold k$ a field. $R$ is a Koszul ring
iff Tor$_i^R(\bold k,\bold k)$
has pure degree $i$ for all $i$.
\endproclaim

We recall now a cohomological interpretation, due to
Lazarsfeld, of the Koszul
property for a coordinate ring  $R(L)$. Let $L$ be a
globally  generated line
bundle on a variety
$X$. We will denote $M^{0,L} := L$ and $M^{1,L} :=M_L \otimes
L = M_{M^{0,L}} \otimes L$. If $M^{1,L}$ is globally generated,
we denote 
$M^{2,L}  :=M_{M^{1,L}} \otimes L$.  We repeat the process
and define inductively $M^{h,L}  :=M_{M^{h-1,L}} \otimes L$,
if $M^{ h-1,L}$ is
globally generated. Now we are ready to state 
the following slightly modified version of \Pareschione \ns, Lemma 1:

\proclaim {\KLaz}
Let $X$ be a projective variety over an algebraic closed field
$\bold k$. Let $L$ be an ample and  base-point-free line
bundle on
$X$. Then $R(L)$ is Koszul iff $M^{h,L}$ exists, is globally
generated and
$H^0(M^{h,L}) \otimes H^0(L^{\otimes s+1}) \to H^0(M^{h,L} \otimes 
L^{\otimes s+1})$ is surjective for all $h \geq 0$, $s \geq 0$. If, in
addition, $H^1(L^{\otimes s+1})=0$ for every $s \geq 0$, then $R(L)$ is
Koszul iff 
\text{H$^1(M^{(h),L}
\otimes L^{\otimes s}) =0$} for every $h \geq 0$ and every $s \geq
0$. 
\endproclaim

The proof of \EKoszul will follow the same strategy of Section $2$, i.e.,
we will translate the problem in terms of a question about vector
bundles over a suitable curve $C$ of $X$. For that purpose we need
now a way to relate $M^{(h),L}$ to $M^{(h),L\otimes \Cal O_C}$.  We
carry this out link by link:

\proclaim{\defrK} Let $X$ be a variety, let $L$ be a line bundle on
$X$ and let
$\frak b$ be a (smooth) effective divisor on $X$. Assume that $M^{h',L}$
is defined for all $h \geq h' \geq 0$ (i.e., inductively, $M^{h'-1,L}$ is
defined and globally generated). We then define, for all $0 \leq h' \leq
h$, 
$M^{h',L}_{h',\frak b} = M^{h',L}
\otimes \Cal O_{\frak b}$. Then $M^{h',L}_{h',\frak b}$ is globally
generated and we define $M^{h'+1,L}_{h',\frak b} = M_{M^{h',L}_{h',\frak
b}}\otimes L$. If $M^{h'+1,L}_{h',\frak b}$ is again globally generated
we define $M^{h'+2,L}_{h',\frak b} = M_{M^{h'+1,L}_{h',\frak
b}}\otimes L$ and so on.
\endproclaim

\proclaim{\Ksequence}
Let $X$ be a variety, let $\frak b$ be a (smooth) effective divisor
on $X$  and let 
$B = \Cal O(\frak b)$.  Let $L$ be a base-point-free line bundle on
$X$
 such that $M^{h',L}$ is globally generated and $H^1(M^{h',L}
\otimes B^*)=0$ for all
$0
\leq h'
\leq h-1$,
$H^1(L
\otimes
\Cal O_\frak b)=0$, and
$L
\otimes \Cal O_\frak b$ is Koszul.
Then,
\roster
\item"(1)"  $M^{h,L}_{h',\frak b}$ is globally generated for all $0 \leq h'
\leq h$.
\item"(2)" $H^1(M^{h,L}_{h',\frak b})=0$ for all $0 \leq h'
\leq h$.
\item "(3)" $0 \to H^0(M^{h'-1,L} \otimes B^* )\otimes
M^{h-h',L}_{0,\frak b} \to M^{h,L}_{h',\frak b} \to M^{h,L}_{h'-1,\frak b}
\to 0$, for all $1 \leq h' \leq h$.
\endroster
\endproclaim

\noindent {\it Proof.} The proof is by induction on $h$. If $h=0$, the
result is part of the hypotheses. If $h =1$, the exact sequence in
$(3)$ is
\Splemmaone
when we set  $F= L$ and twisted by $L$. Let us write $L_\frak b = L
\otimes
\Cal O_\frak b$. Since
$H^1(L_\frak b)=0$ and $L_\frak b$ is Koszul, 
$H^1(M^{1,L}_{0,\frak b})=0$, therefore using
$(3)$ we obtain indeed that $H^1(M^{1,L}_{1,\frak b})=0$. The
bundle $M^{1,L}_{0, \frak b}$ is globally generated because
$L_\frak b$ is Koszul. Finally the fact that
$M^{1,L}_{1,\frak b}$ is globally generated follows again from
$(3)$: we have the following exact commutative diagram
$$\matrix
  H ^0(L \otimes B^*) \otimes H^0(L _\frak b) \otimes  \Cal
O _\frak b &
\hookrightarrow & H^0(M_L \otimes L _\frak b) \otimes \Cal O _\frak b
&\twoheadrightarrow & H^0( M_{L _\frak b} \otimes
L_\frak b) \otimes
\Cal O _\frak b \\ \downarrow && \downarrow && \downarrow
\\  H ^0(L \otimes B^*) \otimes L_\frak b &
\hookrightarrow & M_L \otimes L_\frak b
&\twoheadrightarrow &  M_{L_\frak b} \otimes L_\frak b \\
\endmatrix$$
in which the vertical side arrows are surjective because $L_\frak b$ 
and $ M_{L_\frak b} \otimes
L_\frak b = M^{1,L}_{0,\frak b}$ are both globally generated.
Let us now assume the result to be true for $h-1$ and prove it for
$h$. We again prove  $(3)$ first.  
If $h=h'$, again
$(3)$ is nothing but \Splemmaone \ns,  setting 
$F= M^{h-1,L}$ (which we know by induction hypothesis to be
globally generated) and twisted by $L$. If $h > h'$, by induction on
$h$ we have the sequence
$$0 \to \text H ^0(M^{h'-1,L} \otimes B^*) \otimes
M^{h-h'-1,L}_{0,\frak b} 
\to  M_{h',\frak b}^{h-1,L}
\to   M_{h'-1,\frak b}^{h-1,L} \to  0 \ .$$ Call $V=  H ^0(M^{h'-1,L}
\otimes B^*)$. Taking global sections,   we build this exact
commutative diagram:
$$\matrix
V  \otimes
H^0(M^{h-h'-1,L}_{0,\frak b}) \otimes \Cal O_\frak b  &
\hookrightarrow & H^0(M_{h', \frak b}^{h-1,L})  \otimes
\Cal O _\frak b & \twoheadrightarrow & H^0( M_{h'-1,\frak b}^{h-1,L})
\otimes
\Cal O _\frak b\\ 
\downarrow && \downarrow &&
\downarrow
\\ V \otimes 
M^{h-h'-1,L}_{0,\frak b} &
\hookrightarrow & M_{h', \frak b}^{h-1,L}
&\twoheadrightarrow & M_{h'-1,\frak b}^{h-1,L} \\
\endmatrix$$
The top horizontal sequence is exact at the right because
$H^1(M_{0,\frak b}^{h-h'-1,L})=0$, by induction hypothesis. The
vertical arrows are surjective because the vector bundles involved
are globally generated by induction hypothesis on
$h$. The short exact sequence of kernels is then, after
tensoring by $L_\frak b$, the sequence wanted
for
$(3)$. To prove  $(2)$, we use induction on $h'$. If $h' =0$
both $(1)$ and $(2)$ follow from the fact that $L_\frak b$ is Koszul
and $H^1(L_\frak b)=0$. Now assume that $(1)$ and $(2)$ hold for
$h'-1$.
Condition $(2)$ is a straight forward consequence of already proven
$(3)$ and induction hypothesis on both $h$ and $h'$. 
For $(1)$ we use induction on both $h$ and $h'$
and  $(3)$ just proven.  
If $h=0$ the surjectivity just follows from the fact that 
$L_\frak b$ is Koszul, hence normally generated. If $h'=0$
the surjectivity just follows from the fact that $L_\frak b$
is Koszul. Assume now that the claim holds for $h'-1$. The
surjectivity of the map for $h'$ follows then by chasing the
commutative diagram of multiplication maps, built upon $(3)$, having
in account the vanishing of $H^1(M^{h-h',L}_{0,\frak b})$, which
follows from $(2)$, and the surjectivity of the vertical side maps,
which follows from induction hypothesis on $h$ and $h'$. Then the
fact that $L_\frak b$ is ample implies the global generation of
$M^{h,L}_{h',\frak b}$ as wished. $\square$

\vskip .2 truecm
We are now ready to prove the main theorem of this section:

\proclaim{\EKoszul}
Let $X$ be an Enriques surface. Let $B_1$ and $B_2$ be ample and
base-point-free line bundles, such that $B_1 \cdot B_2 \geq 6$. If $L =
B_1
\otimes B_2$, then $R(L)$ is Koszul.
\endproclaim

\noindent {\it Proof.}
According to \KLaz we need to show that $M^{h,L}$ is globally
generated and that 
$$H^0(M^{h,L})\otimes H^0(L^{\otimes s}) @>\alpha>> H^0(M^{h,L} 
\otimes L^{\otimes s})$$
surjects for all $h \geq 0$ and $s \geq 1$. To better carry out the
argument, is convenient to also prove $H^1(M^{h,L}\otimes B_1^*) =
H^1(M^{h,L} \otimes B_2^*) = 0$.  The proof is by induction on
$h$. If
$h=0$ the result is the projective normality of
$L = B_1
\otimes B_2$, which follows from \Enriquesng \ns, and Kodaira
vanishing.  Now assume the result for $h-1$. Since $L$ is ample, the
surjectivity of $\alpha$ implies the global generation of
$M^{h,L}$, hence we can assume that $M^{h',L}$ is globally
generated for all $0 \leq h' < h$ and we need only to prove
that $\alpha$ surjects and that $H^1(M^{h,L}\otimes B_1^*) =
H^1(M^{h,L} \otimes B_2^*) = 0$. We start proving the
former and in the course of the proof we will also obtain the
desired vanishings. According to
\sfreeobs we are done if we prove that certain collection of
multiplication maps surject. We prove the surjectivity of the first of
them, which is  
$$H^0(M^{h,L})\otimes H^0(B_1) @>\beta>> H^0(M^{h,L} \otimes B_1) 
\ .$$
The argument to prove the surjectivity of the rest is analogous. We 
prove it using again induction on $h$. We proved the statement for $h
=0$ in the course of proving the projective normality of $L$ in
\Enriquesng \ns. Assume the statement to be true for $h-1$ (we may
also assume the surjectivity of the map
$\beta$ for $h-1$ if we substitute in the formula $B_1$ by $B_2$,
since the roles of $B_1$ and $B_2$ are interchangeable.
Consider the sequence
$$\displaylines{
H^0(M^{h-1,L})\otimes H^0(B_2) @>\gamma>> H^0(M^{h-1,L}\otimes 
B_2)\cr
\to H^1(M^{h,L} \otimes B_1^*) \to H^1(M^{h-1,L})\otimes 
H^0(B_2)\ .}$$
The multiplication map $\gamma$ is surjective by induction 
hypothesis. The group $H^1(M^{h-1,L})$ vanishes also by induction
hypothesis, therefore $ H^1(M^{h,L} \otimes B_1^*)=0$. On the other
hand $H^1(\Cal O_X)=0$, so in order to see the surjectivity of
$\beta$ it is enough to check the surjectivity of 
$$ H^0(M^{h,L}_{h,\frak b_1})\otimes H^0(B_1 \otimes 
\Cal O_{\frak b_1}) @>\delta>> H^0(M^{h,L}_{ h,\frak b_1}\otimes B_1)\
,$$
where $\frak b_1$ is a smooth irreducible curve in $|B_1|$. To see the
surjectivity of
$\delta$ we will use
\Ksequence inductively on $h'$. More precisely we want to prove that
$$ H^0(M^{h,L}_{h',\frak b_1})\otimes H^0(B_1 \otimes \Cal 
O_{\frak b_1}) \to H^0(M^{h,L}_{ h,\frak b_1}\otimes B_1) $$
surjects for all $0\leq h'\leq h$. If $h'=0$, $M^{h,L}_{0,\frak b_1}$ is 
semistable with slope strictly bigger than $2g+2$ by \slopelemma \ns,
hence by
\Butres the multiplication map in question is surjective. Now assume
the statement to be true for $h'-1$. We take global sections in the
sequence in part $(3)$ of the statement of \Ksequence and tensor with
$U = H^0(B_1 \otimes \Cal O_{\frak b_1})$ to obtain the following exact 
commutative diagram,
$$\matrix
W \otimes  H^0(M^{h-h',L}_{0, \frak b_1}) 
\otimes U&
\hookrightarrow & H^0(M^{h,L}_{h',\frak b_1}) \otimes U & \twoheadrightarrow &
H^0(M^{h,L}_{h'-1,\frak b_1}) \otimes U
\cr 
\downarrow &&\downarrow&&\downarrow \cr
W \otimes  H^0(M^{h-h',L}_{0,\frak b_1} \otimes B_1) &
\hookrightarrow & H^0(M^{h,L}_{h',\frak b_1} \otimes B_1) &
\twoheadrightarrow & H^0(M^{h,L}_{h'-1,\frak b_1} \otimes B_1) \ ,\cr
\endmatrix $$
where $W = H^0(M^{h'-1,L} \otimes B_1)$. The surjectivity of the left
hand side vertical map and the exactness  at the right of the top
horizontal sequence follow both from
\Butres and \slopelemma \ns. The surjectivity of the right hand side
vertical map follows by the induction hypothesis on $h'$. $\square$

\heading  4. Abelian and bielliptic surfaces \endheading 
In this section we deal with the remaining classes of
surfaces with Kodaira dimension $0$, namely, those
with nonzero  irregularity. For
the techniques employed we return 
to those used in the
arguments of Section 1.  The main theorem we 
will prove is

\proclaim{\ab}
Let $X$ be an Abelian or a bielliptic surface. 
Let $B$ be
an ample and base-point-free line bundle with 
$B^2 \geq 5$ and let $N$ be a numerically
trivial line bundle on $X$. Let
$L_1 \equiv B^{\otimes l_1+1}$ and $L_2 \equiv 
B^{\otimes
l_2+1}$. If
$l_1, l_2
\geq p \geq 1$, then
$$H^1(M_{L_1}^{\otimes p+1} \otimes L_2) =
H^1(M_{L_1} \otimes L_2) =0 
\ .$$
In particular, if  
$n
\geq p \geq 1$, then $B^{n+1} \otimes N$
satisfies the property
$N_p$.

\endproclaim

Before we prove \ab we need the following 

\proclaim{\ablemma}
Let $X$ be a surface with $\kappa =0$. Let $B$ be an 
ample and
base-point-free line bundle. Let $L_1=B^1_1 \otimes 
\dots \otimes B^1_{l_1}$,
where $B^1_i \equiv
B$ are base-point-free line bundles and $l_1 \geq 1$ 
and
$L_2 =B^2_1
\otimes
\dots
\otimes B^2_{l_2}$, where
$B^2_j
\equiv B$ and $l_2 \geq 1$. If either
\roster 
\item "(1)" 
 $l_1$ or $l_2$ are greater than or equal to $3$ or,
\item"(2)"  $l_1=2$ , $l_2=1 \text{ or } 2$ and $H^2(L_1 \otimes
(B^2_1)^{-2})=0$ or,
\item "(3)" $X$ is Abelian or bielliptic surface, $B^2 \geq 5$, and $l_1 =
l_2 =2$,
\endroster
then
$H^1(M_{L_1}\otimes L_2)=0$.
\endproclaim

\noindent {\it Proof.}
In cases (1) and (2) the result follows from 
iteratively 
applying \CM using \sfreeobs and Kodaira vanishing.
In  case (3),  let us
write
$L_1$ as $B^{\otimes 2} \otimes E_1$ with $E_1\equiv
0$. We can find  $E\equiv 0$ with $E^{\otimes 2} 
\neq E_1 \otimes K^*$,
because not all elements in Pic$^0(X)$ (the group of all 
numerically
trivial line bundles up to linear equivalence) have
order
$2$. Then, by \numfree \ns, we can  assume that
$B_1^2 = B
\otimes
E$.  Then $H^2(L_1 \otimes (B_1^2)^{-2}) = H^0(K\otimes
E^{\otimes 2}
\otimes E_1^*)=0$, which follows from our choice of $E$,
and we are in case $(2)$.
$\square$

\vskip .2 truecm

\noindent (4.3) {\it Proof of \ab \ns.} 
The vanishing of $H^1(M_{L_1}\otimes L_2)$ is 
a consequence of \ablemma \ns. The proof of the
vanishings of
$H^1(M_{L_1}^{\otimes p+1}\otimes L_2)$ is by
induction. As usual the key step is the first:
$p=1$. We need to prove that $H^1(M_{R_1}^{\otimes 2}
\otimes R_2)=0$ if $R_1 \equiv B^{\otimes r_1}$ and
$R_2 \equiv B^{\otimes r_2}$ and $r_1, r_2 \geq 2$.
Using the sequence \sequence we obtain
$$\displaylines{H^0(M_{R_1}\otimes R_2) \otimes H^0(R_1) 
@>\alpha >>
H^0(M_{R_1}\otimes R_1\otimes R_2) \cr
 \to H^1(M_{R_1}^{\otimes 2}
\otimes R_2) \to H^1(M_{R_1}\otimes R_2) \otimes H^0(R_1) 
\ .}$$
The  group
$H^1(M_{R_1}\otimes R_2)$ vanishes by \ablemma \ns.
Therefore the sought vanishing is equivalent to the
surjectivity of 
$\alpha$. If $r_1 \geq 3$, let $B_1^1=B$. If $r_1 =
2$, let
$R_1 = B^{\otimes 2}
\otimes E_1$. Analogously, if $r_2 = 2$, let
$R_2 = B^{\otimes 2} \otimes E_2$. We may now assume
if $r_1 \geq 2$, by
\numfree \ns, that
$B_1^1 = B \otimes E$ with $E \equiv 0$ but $E^{\otimes 2}
\neq K^*\otimes E_2$ and, if in addition $r_2=2$,
that 
$E^{\otimes 2} \neq K \otimes E_2^{\otimes 2} \otimes
E_1^*$ also. We can always find such an $E$ if not all
elements in Pic$^0(X)$ have order $2$, $4$ or $6$.
This is the case for Abelian and bielliptic surfaces,
which possess numerically trivial line bundles of
infinite order. Then, to  see that
$\alpha$ surjects, by
\CM and
\sfreeobs it suffices to  check that
$$\displaylines{ 
\abfive 
\hfil H^1(M_{R_1}\otimes R_2 \otimes
(B_1^1)^*) =  0 \hfil \cr
\absix \hfil H^2(M_{R_1} \otimes R_2 \otimes
(B_1^1)^{-2}) =0 \hfil \cr
\abseven \ \ \ \ \hfil H^1(M_{R_1}\otimes R'_2 \otimes
B^{\otimes \gamma}) =  0 \text{
for all } 0 \leq \gamma  \text{ and } R'_2
\equiv R_2 \hfil
\cr
 \ \abeight \hfil H^2(M_{R_1} \otimes R"_2 \otimes
B^{\otimes \gamma -1} ) =0 \text{
for all } 0 \leq \gamma  \text{ and } R"_2
\equiv R_2 \ .\hfil
\cr}$$ The vanishing \abeight follows from \sequence  and
Kodaira vanishing Theorem.
The vanishing \absix follows from \sequence \ns,
Kodaira vanishing Theorem and the way in which we have chosen
$E$.  The vanishing in \abseven follows from \ablemma 
\ns. Finally,  
\abfive follows from
\ablemma once we see that if $r_1=r_2=2$, $H^2(R_1 \otimes
R_2^{-2}
\otimes (B_1^1)^{\otimes 2})=H^2(E_1\otimes E^{\otimes
2} \otimes E_2^{-2})=0$, which follows from the way in
which we  have chosen
$E$.

Assume the result true for $p-1$ and $p >1$. We have
the following sequence:
$$\displaylines{ H^0(M^{\otimes p}_{R_1} \otimes  R_2)
\otimes H^0(R_1) @>\beta >> H^0(M^{\otimes p}_{R_1}
\otimes 
 R_1 \otimes R_2) \cr
\to H^1(M_{R_1}^{\otimes p+1} \otimes R_2) \to
H^1(M^{\otimes p}_{R_1} \otimes R_2) \otimes H^0(R_1) \
.}$$
The last term is zero by induction hypothesis, so the
desired vanishing is equivalent to the surjectivity of
$\beta$. This follows from \sfreeobs and \CM \ns. In
fact, the required vanishings follow by induction,
Kodaira vanishing Theorem and \ablemma \ns. 

For the last conclusion of the theorem, note that 
$$H^1(M_{L}^{\otimes p'+1} \otimes L^{\otimes s}) = 0 \text{ for all } p
\geq p'
\geq 0 \text{ and all } s \geq 1\ .$$
Then,  by \GLlemma \ns, 
$L$ satisfies
property $N_p$.
 
\vskip .2 truecm
Either as  straight forward consequence of \ab or from
the same ideas we have been using we obtain results 
for adjoint linear series:

\proclaim{\abcor} Let $X$ be an Abelian or  bielliptic
surface. Let $B$ be an ample and base-point-free line
bundle such that $B^2 \geq 5$. Then $K \otimes B^{\otimes
n}$  satisfies property $N_p$ if $n \geq p+1$, $p \geq 1$.
\endproclaim
 
\abcor  implies  Mukai's conjecture  for
Abelian and bielliptic surfaces and $p=0$, and for $p=1$ 
(in the latter
case, our result improves Mukai's bound):

\proclaim{\abMukai} Let $X$ be an Abelian or a bielliptic 
surface. Let
$A$ be an ample line bundle and $L = K_X \otimes A^{\otimes 
n}$. If $n
\geq 2p+2$ and $p \geq 1$, then $L$ satisfies property $N_p$. In
particular, if $n \geq 4$, $L$ satisfies property $N_1$.
\endproclaim

\noindent {\it Proof.}
$A^{\otimes 2}$ is base-point-free by \twoamplearefree (for Abelian surfaces this
also follows from Lefschetz's Theorem) and
 since
$K
\equiv 0$, $A^2
\geq 2$ and  $(A^{\otimes 2})^2 \geq 8$. Then, if $n$ is even the
result is a straight forward consequence of \abcor \ns. If $n$ is
odd the situation is the same as that of \EnriquesMukaiNp and we
proceed analogously. 
$\square$
\vskip .2 truecm

\noindent {\bf \abKempf \ns.} If $X$ is an Abelian
surface the above   result was proven by Kempf (cf.
\Kempf
\ns). 
However, the results proven in this chapter are  more general: for
instance, since on an Abelian surface a polarization of type 
$(1,3)$ is
base-point-free, \abcor implies that a line bundle of type $(p+1,
3p+3)$ satisfies property $N_p$. This fact does not follow from
Kempf's result.
\vskip .2 truecm
To end this section we carry out a study analogous to the one 
realized for Enriques surfaces in Section $3$: the
following theorem proves in particular  that
the  line bundles satisfying property
$N_1$ according to \ab  have also a Koszul coordinate
ring.

\proclaim{\abKoszul}
Let $X$ be an Abelian or a bielliptic surface. 
Let  $B_1$ and $B_2$ be
numerically equivalent ample and base-point-free line bundles with 
self-intersection bigger than or equal to $5$. If $L=B_1 \otimes 
B_2$,
then $R(L)$ is Koszul. In particular $L$ satisfies
property $N_1$.
\endproclaim

In order to prove the theorem we use the following result which is 
basically a reformulation of \GPone \ns, Theorem 5.4 for the case of
surfaces with $\kappa =0$:

\proclaim{\abKlemma}
Let $X$ be a surface with $\kappa = 0$, let $B_1$ and $B_2$ be two
ample and base-point-free line bundles. If 
$ H ^2(B_1
\otimes B_2 ^*) = H ^2(B_2 \otimes B_1 ^*)=0$,  
then the following properties are satisfied for all $h \geq 0$:
\vskip .3 cm
$1)$ \ $M^{h,L}$ is globally generated 
\vskip .05 cm
$2)$ \ H$^1(M^{h,L} \otimes B_1 ^{\otimes b_1} \otimes B_2
^{\otimes b_2})=0$ for all $ b_1, b_2 \geq 0$
\vskip .05 cm
$3)$ \ H$^1(M^{h,L} \otimes B_j
^*)=0$ where $j =1,2$ 
\vskip .05 cm
$4)$ \ H$^1(M^{h,L}\otimes B_i \otimes B_j
^*)=0$ where $i =1,2$ and $j=2,1$
\vskip .05 cm
$5)$ \ H$^1(M^{h,L}\otimes B_i ^{\otimes 2} \otimes B_j
^*)=0$ where $i =1,2$ and $j=2,1$
\vskip .3 cm
In particular H$^1(M^{(h),L} \otimes L
^{\otimes s})=0$ for all $h,s \geq 0$, and $R(L)$ is a 
Koszul
$k$-algebra.
\endproclaim

\noindent {\it Proof.}
For the proof of the lemma we refer to \GPone \ns. Since now $B_1$ 
and $B_2$ are ample and $K_X \equiv 0$, we obtain all the 
vanishings
of the groups $H^1(B_1^{\otimes a} \otimes B_2^{\otimes
b})$  when $a,
b \geq 0$ and $a+b\geq 1$ needed in the proof, by Kodaira 
vanishing
Theorem, making therefore unnecessary to assume the 
vanishings of
$H^1(B_1)$, $H^1(B_2)$ and $H^2(\Cal O_X)$.  $\square $
\vskip .2 truecm
\noindent 4.9. {\it Proof of \abKoszul \ns.}
The result follows from \abKlemma \ns. If 
$X$ is bielliptic, the only problem we might have is if
$B_1=B_2 \otimes K_X$ or if $B_2=B_1 \otimes K_X$. In the
former case, choose a line bundle $E \in
\text{Pic}^0(X)$ such that $E ^{-2} \neq 0$ and
$K_X^{\otimes 2}  \otimes E ^{\otimes 2} \neq 0$. In
the latter case choose $E$ such that $E
^{\otimes 2}\neq 0$ and $K_X^{\otimes 2}  \otimes E
^{-2} \neq 0$. Let then $B_1'=B_1 \otimes E$ and
$B_2'=B_2 \otimes E^*$. The desired result follows
if we apply \abKlemma to $B_1'$ and $B_2'$ instead.

If $X$ is an Abelian surface, $K_X$ is trivial, so  the only problem applying \abKlemma would appear when $B_1=B_2$. This is solved analogously considering $B_1'=B_1 \otimes E$ and $B_2'=B_2 
\otimes E ^*$, where now $E$ is taken to have nontrivial
square.
$\square$

\heading 5. Surfaces of positive Kodaira dimension \endheading

In this section we focus on the study 
of  adjoint linear series of surfaces of positive Kodaira
dimension. We find sufficient conditions for the normal generation
and the normal presentation of the adjoint linear series and of the
powers of an ample and base-point-free line bundle. For the latter
case we also generalize the results to higher syzygies. One can look
upon these results as an analogue for projective normality and higher
syzygies of the results of Kawamata and Shokurov (see \Kaw and \Sho
\ns) for base-point-freeness and effectiveness, viewed  in the
special context of algebraic surfaces. They deal with nef bundles $L$
for which $L \otimes K^*$ is nef and big and  conclude the freeness
and effectiveness of multiples of $L$. We start with an ample and
base-point-free bundle $B$ which satisfies certain
inequalities (see Theorems 5.1 and 5.8) which are immediate if one
assumes that $B \otimes K^*_S$ is nef and big and go on to prove
projective normality and higher syzygy results for powers of $B$ and
for adjunction bundles associated to $B$. 

We obtain two interesting
consequences from our study. The first is finding   sufficient
conditions for projective normality and quadratic generation of
pluricanonical embeddings of surfaces of general type (\fourK
\ns, Remark 5.7 and \fiveK
\ns).  Bombieri asked in \Bombieri whether $|K_S^{\otimes 5}|$
maps $S$ as a projectively normal variety.
This question 
was answered affirmatively by Ciliberto in
\Ciliberto (under basically the same assumptions of
\twoK below). Thus the above mentioned corollaries recover,  and in the
case of regular surfaces, improve Ciliberto's result.   
The second consequence
is an effective bound along the lines of Mukai's
conjecture using a result by Fern\'andez del Busto. 
In the case of pluricanonical models of regular surfaces of general
type we further our study to higher syzygies.   Ein and Lazarsfeld's
results in \EL together with Del Busto's give effective bounds
(slightly weaker than ours) along the lines of the Mukai's
conjecture, but for regular surfaces, the bounds we obtain are 
better. We also obtain as a corollary of Theorem 5.1 (4)  
effective bounds for property $N_p$ for the
multiples of ample line bundles on regular surfaces.

\proclaim{Theorem 5.1 }
Let $S$ be a regular surface of positive Kodaira dimension and
$p_g
\geq 4$. Let $B$ be an ample and base-point-free line bundle such
that
$H^1(B)=0$. Let $L=K_S
\otimes B^{\otimes n}$ and $L'=K_S
\otimes B^{\otimes l}$. Let $N=B^{\otimes m}$ and $N'=B^{\otimes
k}$.
\roster
\item"(1)"If $\kappa(S) =1$ and $B^2>K_S \cdot B$, and if
$n,l \geq 2$, then
$H^1(M_L\otimes L')=H^1(M_L^{\otimes 2}\otimes L')=0$. In
particular, $K_S
\otimes B^{\otimes n}$ satisfies property $N_1$, for all $n \geq
2$.
\item"(2)" If $\kappa(S)=2$ and $B^2 \geq K_S \cdot B$, and if
$n,l \geq 2$, then
$H^1(M_L\otimes L')=0$. In particular, 
 $K_S
\otimes B^{\otimes n}$ satisfies property $N_0$, for all $n \geq
2$.
\item"(3)" If $\kappa(S)=2$ and $B^2 \geq 2 K_S \cdot B$, and if 
$n,l
\geq 2$, then $H^1(M_L\otimes L')=H^1(M_L^{\otimes 2}\otimes
L')=0$; and if $m,k \geq 2$, then $H^1(M_N\otimes
N')=H^1(M_N^{\otimes 2}\otimes N')=0$. In particular
$K_S
\otimes B^{\otimes n}$ and $B^{\otimes m}$ satisfy property
$N_1$, for all
$n,m
\geq 2$.
\item"(4)" If $\kappa(S)=2$ and $B^2 \geq 2 K_S \cdot B$,  and if
$m,k \geq p+1, p \geq 1$, then $H^1(M_N^{\otimes p+1} \otimes
N')=0$. In particular if $p \geq 1$, 
$B^{\otimes m}$ satisfy property
$N_p$, for all
$m
\geq p+1$.
\endroster 
\endproclaim

To prove Theorem 5.1 we will need these two lemmas:
\proclaim{Lemma 5.2}
Let $S$ be a surface 
and $B$ an ample and base-point-free line bundle with
$H^1(B)=0$ and $B^2 > B \cdot K_S$ if $\kappa(S) \leq 1$, and $B^2 \geq B \cdot
K_S$ if $\kappa(S)=2$. Then
$H^1(B^{\otimes m})=0$ for all
$m \geq 1$.
\endproclaim

\noindent{\it Proof}.  Let $C$ smooth curve in $|B|$. Since deg
$(B^{\otimes m} \otimes \Cal O_C) > 2g(C)-2$ when $m \geq 3$, we only
have to prove $H^1(B^{\otimes 2})=0$. If $B^2 >B \cdot K_S$ or
$B^{\otimes 2} \otimes \Cal O_C \neq K_C$, then
$H^1(B^{\otimes 2}\otimes \Cal O_C)=0$, hence 
$H^1(B^{\otimes 2})=0$ because $H^1(B)=0$. If $B^{\otimes 2} \otimes
\Cal O_C = K_C$, then $B \otimes \Cal O_C = K_S \otimes \Cal O_C$.
Consider the sequence 
$$0 \to H^0(K_S^*) \to H^0(B \otimes K_S^*) \to H^0(B \otimes K_S^*
\otimes \Cal O_C) \to H^1(K_S^*) \ .$$
Since in this case $S$ is a surface of general type,
$H^0(K_S^*)=H^1(K_S^*)=0$, therefore $B \otimes K_S^*$ is
effective and since $B$ is ample, it must be $B \otimes K_S^*=\Cal
O_S$. Hence $H^1(B^{\otimes 2})=H^1(K_S^{\otimes 2})=0$. $\square$

\proclaim{Lemma 5.3}
Let $S$ be an algebraic surface with
nonnegative Kodaira dimension and let $B$ be an ample line
bundle. Let $m \geq 1$. If $B^2 \geq mK_S \cdot B$, then $K_S \cdot B
\geq mK_S^2$.
\endproclaim

\noindent {\it Proof}.
We assume the contrary, i.e., that $K_S \cdot B < mK_S^2$,  
and get a contradiction. Let $L=B
\otimes K_S^{ -m}$. We have that $L^2 >0$. By Riemann-Roch 
$$h^0(L^{\otimes n}) \geq \frac{n^2L^2-nK_S\cdot L}{2} + \chi(\Cal O_S)
- h^0(K_S \otimes L^{ -n}) \ .$$ If $B^2 > mK_S \cdot B$, $(K_S
\otimes L^{\otimes -n})\cdot B <0$, for $n$ large enough, and since
$B$ is ample, 
$K_S
\otimes L^{\otimes -n}$ is not effective, so finally $L^{\otimes n}$
is effective for $n$ large enough. But in that case $nK_S \cdot L \geq
0$, because $K_S$ is nef, contradicting our assumption.

Now if $B^2 = mK_S \cdot B$, we have that $L^2 >0$, $B^2>0$ (because
$B$ is ample), and $L \cdot B =0$, but this is impossible by the Hodge
index theorem.
$\square$

\vskip .2 truecm
\noindent (5.4) {\it Proof of Theorem 5.1}.
We start proving that $H^1(M_L\otimes L')=0$ if $B$ satisfies the conditions
of (1) and (2). 
By \sfreeobs it suffices to show that
$$\displaylines {H^0(K_S \otimes B^{\otimes n}) \otimes H^0(B)
\longrightarrow H^0(K_S \otimes B^{\otimes n+1}) \cr
H^0(K_S \otimes B^{\otimes m}) \otimes H^0(K_S \otimes B)
\longrightarrow H^0(K_S \otimes B^{\otimes m+1}), \text{ for all
}n \geq 2, m \geq 3 }$$
surject. We want to show now that $K_S \otimes B$ is base-point-free. Let $C
\in |B|$ smooth. From the
arguments of the proof of Lemma 5.2, we see that $h^2(B) \leq 1$ with
equality if and only if $B=K_S$ and $S$ of general type.
Since $S$ is regular we have
$$0 \to H^0(\Cal O_S) \to H^0(B) \to H^0(B \otimes \Cal O_C) \to 0 \ .$$
Thus it follows from this and from Riemann-Roch that $h^0(B \otimes \Cal O_C)
\geq p_g -1$ (again, with equality if and only if $B=K_S$ and $S$ of general
type). Now by Clifford's bound, $B^2 \geq 6$ except if $B=K_S$ and $S$ of
general type. It is well known \text{(cf. \twoK \ns)} that
$K_S^{\otimes 2}$ is base-point-free under the hypothesis of the
theorem. Now, if
$B^2 \geq 6$,
$K_S \otimes B$ is base-point-free by \numfree \ns.
Take now  
$C' \in |K_S \otimes B|$, also a smooth curve.
Since
$H^1(\Cal O_S)=0$ and by Lemma 5.2 and Kodaira vanishing theorem, we
may apply
\robs
and conclude that it is  enough to check that
$$\displaylines {H^0(K_S \otimes B^{\otimes n} \otimes \Cal O_C)
\otimes H^0(B \otimes \Cal O_C) \longrightarrow H^0(K_S \otimes
B^{\otimes n+1} \otimes \Cal O_C)
\cr H^0(K_S \otimes B^{\otimes m} \otimes \Cal O_{C'}) \otimes
H^0(K_S \otimes B \otimes \Cal O_{C'}) \longrightarrow H^0(K_S
\otimes B^{\otimes m+1} \otimes \Cal O_{C'})}$$ surject 
for
all $n \geq 2, m \geq 3$.
This follows from
\Butres or \Pareschilemma \ns. We check this explicitly for the first
family of maps. Let $G= K_S \otimes B^{\otimes n} \otimes \Cal O_C$ and
$G'=B \otimes \Cal O_C$. It is enough to show that deg$G > 2g(C)$ and
that deg$G +
\text{ deg }G' > 4g-2h^1(G')$. For the first inequality, note that $K_S
\otimes B \otimes \Cal O_C= K_C$ and that deg$(B \otimes \Cal O_C) \geq
4$. For the second,  
$$\text{ deg }G + \text{ deg }G' \geq K_S \cdot B + 3B^2 \geq
2(K_S \cdot B +B^2) = 4g(C)-4 \ . \leqno (5.4.1) $$ 
On the other hand
$h^1(G')= p_g-h^0(K_S
\otimes B^*)$ and since $B$ is ample, the inequality $B^2 \geq K_S \cdot B$
implies $h^1(G') \geq p_g-1$. By the bound on $p_g$ we therefore have 
 $4g-4 \geq 4g+2 - 2h^1(G')$. The
reasoning for the second family of maps is similar.

We go on now to prove (3). The proof of the vanishing of 
$H^1(M_N \otimes N')$
uses the same above arguments  and we will not repeat them here.  
We now prove the vanishings of  $H^1(M_L^{\otimes 2} \otimes L')$ and  
$H^1(M_N^{\otimes 2}
\otimes N')$. By  \sfreeobs it is enough to show that the following maps
$$\displaylines{H^0(M_{L} \otimes L') \otimes H^0(B) @>
\alpha_1 >> H^0(M_{L} \otimes L'\otimes B) \cr
H^0(M_{L} \otimes L' \otimes B) \otimes H^0(K_S \otimes B) @>
\alpha_2 >> H^0(M_{L} \otimes L' \otimes K_S \otimes B^{\otimes
2})\cr
H^0(M_{N} \otimes N') \otimes H^0(B) @>
\alpha_3 >> H^0(M_{N} \otimes N'\otimes B)}
$$
surject.
We only sketch in some detail the proof of the surjectivity of
$\alpha_1$, as the proofs for the other two maps are analogous. We
will use \robs and \Splemma \ns. For that we need to check
that
$H^1(M_{L}\otimes L' \otimes B^*)=0$. This  follows from the
surjectivity of
$$\displaylines{H^0(L) \otimes H^0(B) @> \beta_1 >> H^0(L
\otimes B)
\cr
 H^0(L) \otimes H^0(K_S \otimes B) @> \beta_2 >>
H^0(L
\otimes K_S \otimes B) \ .}$$
The surjectivity of $\beta_1$ have already been shown
previously in this proof. We show now the surjectivity of $\beta_2$.
Let $C' \in |K_S \otimes B|$. By the vanishing of $H^1(B)$ from the
hypothesis and Kodaira vanishing the surjectivity of $\beta_2$ follows
from the surjectivity of 
$$H^0(L \otimes \Cal O_{C'}) \otimes H^0(K_S \otimes B \otimes \Cal
O_{C'}) @>\gamma >> H^0(L \otimes K_S \otimes B \otimes \Cal O_{C'})$$
by \robs \ns. We want to apply \Pareschilemma \ns. Note that deg$(L
\otimes
\Cal O_{C'}) +
\text{ deg } (K_S \otimes B \otimes \Cal O_{C'}) \geq
2K_S^2+5K_S\cdot B + 3B^2$ and $4g(C')-4=4K_S^2+6K_S\cdot B+2B^2$.
Since $B^2 \geq 2K_S \cdot B$, by Lemma 5.3 $B^2 \geq 2K^2 + K_S
\cdot B$. Finally $h^1(K_S \otimes B \otimes \Cal
O_{C'}) = p_g \geq 4$, so the inequalities needed to
apply \Pareschilemma are satisfied and the surjectivity of
$\gamma$ follows. 
Returning to the proof of the surjectivity of
$\alpha_1$, we may now apply \robs and \Splemma and therefore it
suffices to check the surjectivity of
$$H^0(M_{L \otimes \Cal O_C} \otimes L' \otimes \Cal O_C)
\otimes H^0(B \otimes \Cal O_C) \longrightarrow H^0(M_{L \otimes
\Cal O_C}
\otimes L'\otimes B \otimes
\Cal O_C) \ ,$$ 
which  follows from \Butres
\ns. Indeed, let $E=B \otimes \Cal O_C$ and $F=M_{L \otimes \Cal O_C}
\otimes L' \otimes \Cal O_C$. We need check that $\mu(F) > 2g(C)$ and that
$\mu(F) > 2g+ \text{ rank }(E)(2g(C)-\mu(E))-2h^1(E)$. For the former
inequality,  $\mu (F) \geq  \mu(M_{L \otimes \Cal O_C}) + K_S \cdot B
+ 2B^2$ and this is bigger than $2g(C)$ since $2g(C)=K_S \cdot B +B^2 +2$
and $B^2 \geq 4$. The latter inequality follows from deg$L \otimes \Cal O_C
>2g(C)$, as $B^2 \geq K_S \cdot B$, $h^1(B \otimes \Cal O_C) \geq 3$ and
$B^2 \geq 4$.

Finally the proof of the vanishing of $H^1(M_L^{\otimes 2} \otimes L')$ if
$\kappa(S)=1$ follows by the same arguments and it is left to the 
reader.

The proof (4) is built upon (3), using induction on $p$ and \CM \ns.
$\square$

\vskip .2 truecm
We want now to apply this theorem to the study of pluricanonical 
models
of surfaces of general type (regular, for the moment; we will
complete the picture when we have at our disposal Theorem 5.8 which
deals with irregular surfaces). The idea is to find a smallest power
of
$K_S$ which is base-point-free, for it would play the role of $B$.  It is
known that under certain mild conditions
$K_S^{\otimes n}$ is base-point-free if $n \geq 2$.
The precise result, which is due to Bombieri, Francia,
Reider and others, can be found in \Catanese
\ns:

\proclaim{\twoK \ns (\Catanese
\ns, Theorem 1.11 (i))} Let $S$ be a surface of general
type. Assume that either
\roster
\item"(1)" $K_S^2 \geq 5$ or
\item"(2)" $K_S^2 \geq 2$ and $p_g \geq 1$, but it
does not happen that $q=p_g=1$ and $K_S^2=3$ or $4$.
\endroster
If $n \geq 2$, then $K_S^{\otimes n}$ is
base-point-free.
\endproclaim
With this we are ready to obtain the following

\proclaim{\fourK}
Let $S$ be a regular surface  of general type with ample canonical bundle
and
$p_{g}\geq 3$. 
Then 
\roster
\item"(1)" $H^1(M_{K_S^{\otimes 4+k}}\otimes K_S^{\otimes
4+l})=0$ for all $k,l \geq 0 $, and 
\item"(2)" 
$H^1(M_{K_S^{\otimes 4+k}}^{\otimes 2}
\otimes  K_S^{\otimes
4+l})=0$ for $k,l=0$, for all $k,l \geq 1$, and for all
$k \geq 0, l \geq 2$.
\endroster 
In particular, if $n \geq 4$, 
 $|K_{S}^{\otimes n}|$ embeds $S$ 
as a
projectively normal variety with homogeneous ideal  
generated by quadratic equations.
\endproclaim

\noindent {\it Proof}. The result is a straight forward consequence of
Theorem 5.1 if
$p_g \geq 4$, setting $B=K_S^{\otimes 2}$. However we can take advantage
of the fact that we are dealing with base-point-free line bundles of
particularly nice shape. If one goes through the steps of the proof of
Theorem 5.1, one sees that one of the places when we use $p_g \geq 4$ is
to prove that $K_S \otimes B$ is base-point-free. In this setting we know
this to be true by Theorem 5.5. The other place where we use the bound on
$p_g$ is when checking the inequalities needed to apply \Pareschilemma
and \Butres \ns. The reader can see that in this particular case $p_g
\geq 3$ suffices for such a purpose.
$\square$

\vskip .2 truecm
\noindent {\bf Remark 5.7.} Some hypothesis of
\fourK can be dropped or relaxed. If we don't require $K_S$ to be
ample, we obtain essentially the same result: the image of $S$ by
$|K_{S}^{\otimes n}|$ is a
projectively normal variety with homogeneous ideal  
generated by quadratic equations. Indeed, note that the ampleness of $B$
was used in Theorem 5.1 to obtain cohomology vanishings and the
base-point-freeness of $K_S \otimes B$. Those are taken care now by
\twoK and Kawamata-Viehweg vanishing theorem. On the other hand we
can relax the hypothesis on $p_g$ to obtain a weaker result, proven by
the same techniques:
\vskip .15 truecm
\noindent{\bf (5.7.1)}{\it 
Let $S$ be a regular surface of general type with
either  $p_g\geq 1$ and 
$K_{S}^{2}\geq 2$ or $K_S^2 \geq 5$. If $n \geq 5$, then
the image of $S$ by the complete linear series 
$|K_{S}^{\otimes n}|$ is a projectively normal variety.}

\vskip .2 truecm

We now complete the picture with the nonzero irregularity case:

\proclaim{Theorem 5.8}
Let $S$ be an irregular surface of positive Kodaira dimension.
Let $B$ be an ample line bundle such that $B^2 \geq 5$ and 
$B'$ is base-point-free and $H^1(B')=0$ for all $B'$ homologous to $B$
(respectively numerically equivalent). Let  $L$ homologous to  $K
\otimes B^{\otimes n}$ (respectively numerically equivalent) and
$L'$
homologous to $K
\otimes B^{\otimes l}$ (respectively numerically equivalent).
\roster  
\item"(1)" If $\kappa(S)=1$ and $B^2 > K_S\cdot B$, and if $n,l \geq 2$,
then 
$H^1(M_L \otimes L')=0$;  if $n,l \geq 3$, then $H^1(M_L^{\otimes 2}
\otimes L')=0$. In particular  
$L$ satisfies property $N_0$ if $n \geq 2$, and  $L$
satisfies property
$N_1$ if $n \geq 3$.
\item"(2)" If $\kappa(S)=2$ and $B^2 \geq 2 K_S \cdot B$, and if 
$n,l \geq 2$,
then 
$H^1(M_L \otimes L')=0$. In particular, $L$
satisfies property $N_0$ if $n \geq 2$;
\item"(3)" if $\kappa(S)=2$ and $B^2 \geq K_S \cdot B$, and if $n,l \geq
3$, then $H^1(M_L
\otimes L')=H^1(M_L^{\otimes 2}
\otimes L')=0$. In particular, $L$ satisfies property
$N_1$ if  for $n \geq 3$.
\endroster
\endproclaim

\noindent{\it Sketch of proof}.
The proof uses  Lemma 5.2, the intersection number inequalities
in our hypothesis and  arguments similar to those in Section 4. We will
outline the argument to show (1).
Assume for simplicity's sake that $L = B^{\otimes n}$. The vanishing of
$H^1(M_L
\otimes L')=0$ is equivalent (because of Kodaira vanishing theorem) to the
surjectivity of 
$$ H^0(L) \otimes H^0(L') \longrightarrow H^0(L \otimes L') \ .$$
Now we break up $L'$ as the tensor product $B_1 \otimes B_2 \otimes \cdots
\otimes B_{l-1} \otimes (K_S \otimes B_l)$ where $B_1 = B \otimes E$ where
$E \in \text{ Pic}^0(S)$ and $E^{\otimes 2} \neq \Cal O$. Clearly, if $L'$
is homologous to $K
\otimes B^{\otimes l}$, $B_i$ is homologous to $B$ and in any case
numerically equivalent, hence by hypothesis and \numfree $B_i$ is
base-point-free and so is $K_S \otimes B_l$. By \sfreeobs it would suffice
to show the surjectivity of several multiplication maps. The first one
would be 
$$H^0(L) \otimes H^0(B_1) \longrightarrow H^0(L \otimes B_1) \ .$$
This follows from \CM, the required vanishings being obtained from Kodaira
vanishing theorem  and, if
$n=2$, from our choice of
$E$. The possible intermediate maps are surjective by \CM and Kodaira. The
last of the maps is surjective \CM \ns, Kodaira vanishing theorem and our
hypothesis. Indeed, the required vanishings of $H^1$ follow from  Lemma 5.2
(here we use the condition $B^2 > B \cdot K_S$), using the divisibility of
Pic$^0$.  We also need the vanishings of 
$H^2(N)$, where $N \equiv B^{\otimes k}
\otimes K_S^*$ and $k
\geq 1$ which certainly occur because $S$ is an elliptic surface.

Now we prove $H^1(M_L^{\otimes 2}
\otimes L')=0$ under the hypothesis in (1). This vanishing is
equivalent to the surjectivity of 
$$H^0(M_L \otimes L') \otimes H^0(L) \longrightarrow H^0(M_L \otimes L
\otimes L')$$
by virtue of the vanishing just shown. If we break $L$ as tensor product
of $K_S$ and base-point-free line bundles homologous or numerically
equivalent to $B$ depending on the case, we can use \sfreeobs and \CM
\ns. The vanishings we need to check have already been shown  in the
first part of the proof.  
$\square$

\proclaim{\fiveK }
Let $S$ be an irregular surface of general type with 
$K_{S}^{2}\geq
5$. If $n \geq 5$, then the image of $S$ by the complete
linear series $|K_S^{\otimes n}|$ is a projectively
normal variety.
\endproclaim

\noindent{\it Proof.}
The observation about ampleness made in Remark 5.7 also applies here.
Having that in account, the  result follows from the arguments of the proof of
Theorem 5.8 if $n$ is odd, and for $n$ even we argue as in Corollary 2.15.
$\square$
\vskip .2 truecm

Another quite interesting consequence of Theorems 5.1 and 5.8, and of a
result  by Fern\'andez del Busto is  the following effective
bound  along the lines of Mukai's
conjecture:
\proclaim{Corollary  5.10}
Let $S$ be an algebraic surface of positive Kodaira dimension,
let $A$ be an ample line bundle and let 
$m = \left [ \frac{(A\cdot(K_S + 4A)+1)^2}
{2A^2} \right ]$. Let
$L = K_S
\otimes A^{\otimes n}$. If
$n
\geq 2m$, then
$L$ satisfies property $N_0$. If $n \geq 3m$, then $L$ satisfies
property $N_1$. If $S$ is a regular surface of general type and $n
\geq 2m$, then $L$ satisfies property $N_1$.
\endproclaim
\noindent{\it Sketch of proof.}
The key observation is the fact that if $k \geq m$ then it follows
from \Busto \ns, Section 2 that 
$A^{\otimes m}$ is base-point-free and $H^1(A^{\otimes m})=0$.
Then we take $A^{\otimes m}$ as the base-point-free line bundle $B$
in Theorems 5.1 and 5.8. One can easily verify that the
numerical conditions in the statements are satisfied.  
Note that Pic$^0(S)$ is divisible, so if $S$ is irregular, Fern\'andez del
Busto's result applies also to $B'$ in the statement of Theorem 5.8.
There is however one hypothesis of Theorem 5.1 which we have
not assumed in this corollary, and which in fact does not occur in
general under the hypothesis of our statement. That is the
assumption of $p_g \geq 4$. This hypothesis was used in the proof
of Theorem 5.1 to check the inequalities needed to apply \Butres or
\Pareschilemma \ns. Under our current hypothesis $B^2$ is much
larger than $K_S \cdot B$ and in any case, large enough to render
the mentioned assumption unnecessary (see for instance (5.4.1)). 
Therefore
the theorem is either a direct consequence of Theorems 5.1 and 5.8 or
follows from slight modifications of the arguments involved in proving those
theorems, for we are in a  situation similar to  Corollary 2.15.
$\square$
\vskip .2 truecm

We focus now on the  study of higher syzygies of
pluricanonical models of surfaces of general type. Recall that one
can obtain a result regarding them from \preg \ns:
\proclaim{\twop+four}
Let $S$ be a surface of general type satisfying the
assumptions of \twoK \ns.
If $n \geq 2p+4$, then the image of $S$ by 
$|K_S^{\otimes n}|$ is projectively normal, its ideal is 
generated
by quadrics and the resolution of its homogeneous 
coordinate
ring is linear until the $p$th stage.
\endproclaim

\noindent {\it Proof.}
The line bundle $K_S^{\otimes 2}$ is base-point-free  by
\twoK \ns.  On the other hand, $K_S^{\otimes 2}$ is
$3$-regular by  the
Kawamata-Viehweg Vanishing
Theorem. Hence  from \preg  the result follows for
$n$ even. 
If $n$ is odd we
argue as in \EnriquesMukaiNp \ns, writing 
$K_S^{\otimes n}$
as $B^{\otimes s-1} \otimes B'$, where 
$B=K_S^{\otimes 2}$ and
$B'=K_S^{\otimes 3}$, which are 
base-point-free by
\twoK \ns. 
$\square$
\vskip .2 truecm

This result can be improved for regular surfaces if we impose the 
hypothesis
of Corollary 5.6: 

\proclaim{\regulargtNp}
Let $S$ be a regular surface of general type with 
$p_{g}\geq 3$.
Let $L= K_S^{\otimes
2p+2+l}$ and
$L'= K_S^{\otimes 2p+2+k}$.    
\roster
\item"(1)" If $p \geq 2$,  $H^{1}(M_{L}^{\otimes
p+1}\otimes L^{\prime })=0$  for all 
$k, l\geq 0$, and
\item"(2)"
if $p=1$, $H^{1}(M_{L}^{\otimes p+1}\otimes
L^{\prime })=$ for $k,l = 0$,  for all
 $k, l \geq 1$,
 and for all $k \geq 0, l \geq 2
$.
\endroster
 
In particular,  if $n \geq 2p+2$ and $p \geq 1$,
then the image of $S$ by 
$|K_S^{\otimes n}|$ is projectively normal, its ideal is 
generated
by quadrics and the resolution of its homogeneous 
coordinate
ring is linear until the $p$th stage. 
\endproclaim

\noindent{\it Proof.}
The proof is by induction on $p$. The statement for  
 $p=1$
is \fourK  $(2)$, having in account the observation on ampleness
made in Remark 5.7. Let us assume the result  to be
true for $p-1$ and 
prove the vanishing for $p$.  Tensoring \sequence
with $M_{L}^{\otimes p}\otimes L^{\prime }$  and taking
global sections yields the following long exact sequence 
$$\displaylines{
H^{0}(L)\otimes 
H^{0}(M_{L}^{\otimes p}\otimes L^{\prime })@>{\eta
}>>H^{0}(M_{L}\otimes L\otimes L^{\prime }) 
\cr
\longrightarrow H^{1}(M_{L}^{\otimes p}\otimes L^{\prime
})\longrightarrow H^{0}(L)\otimes H^{1}
(M_{L}^{\otimes p+1}\otimes
L^{\prime }) \ .}
$$
The last term is zero by induction assumption, thus
the vanishing   is equivalent  to showing
the surjectivity of the multiplication map $\eta$. 
Let $B=K_S^{\otimes 2}$ and $B'=K_S^{\otimes
3}$. By 
\sfreeobs it suffices to
show the surjectivity of several maps:
$$\displaylines{
H^{0}(B)\otimes H^{0}(M_{L}^{\otimes p}\otimes L^{\prime
})@>{\eta ^{\prime }}>>H^{0}(M_{L}\otimes
(B\otimes L^{\prime })), \text{ for all } l,k \geq 0
\text{ and, }\cr
H^{0}(B')\otimes H^{0}(M_{L}^{\otimes p}\otimes L^{\prime
}\otimes B)@>{\eta ''}>>H^{0}(M_{L}\otimes
(B \otimes B'\otimes L^{\prime })), \text{ for all } l,k
\geq 0 \ .}
$$

The surjectivity of $\eta ''$ follows by \CM \ns, the
vanishings required following by induction, from
\sequence and Kawamata-Viehweg.
The surjectivity of $\eta'$ also follows from \CM by
the same reasons, but as in the proof of
\Enriqueshighersyzygies we could alternatively argue
restricting to a smooth curve $C$ in $|B|$. We would have
to eventually use
\Butres and the fact that the needed inequalities hold
follows from adjunction and the assumption on $p_g$,
having in account that $h^1(B \otimes \Cal O_C)=h^2(\Cal
O_S)$.

Finally, the statement on the coordinate ring of the
image of the pluricanonical maps follows from \fourK
\ns, the vanishings just proven and \GLlemma \ns. 
$\square$
\vskip .2 truecm
As a corollary of Theorem 5.1 (4) and Del Busto's result we obtain an
effective bound for a power of an ample line bundle to satisfy
property $N_p$:
\proclaim{Corollary 5.13}
Let $S$ be a regular surface of general type, let $A$ be an ample
line bundle and let 
$m$ as in Corollary 5.10. Let
$L =  A^{\otimes n}$. If
$n
\geq mp+m$, then
$L$ satisfies property $N_p$.
\endproclaim

We state now a result for normal presentation and Koszul
property of adjoint linear series on regular surfaces of general type
with  base-point-free canonical bundle:  

\proclaim{Theorem 5.14} 
Let $S$ be a regular surface of general type with $p_g \geq 4$ and
base-point-free canonical bundle. Let $B$ be an ample and base-point-free line
bundle on
$S$ with $H^1(B)=0$ and let $B^2 \geq B \cdot K_S$. Let 
$L=K_S \otimes B^{\otimes
n}$ and $L'=K_S \otimes B^{\otimes
m}$. If $n,m  \geq 2$, then $H^1(M_L \otimes
L')=H^1(M_L^{\otimes 2} \otimes L')=0$. In particular, if
$n
\geq 2$, then 
$K_S
\otimes B^{\otimes n}$ satisfies  property $N_1$, and in addition, the Koszul
property.
\endproclaim

\noindent {\it Proof}.

{\it Cohomology vanishings}: 
First we check the vanishing
of
$H^1(M_{L}\otimes L^{\prime })$. By Kodaira vanishing theorem, it
suffices to check the surjectivity of 
$$H^0(L^{\prime }) \otimes H^0(L) \to H^0(
L \otimes L^{\prime }) \ .$$
Recall that both $B$ and $K_S$ are base-point-free. By \sfreeobs it suffices
to check the surjectivity of 
$$\displaylines{
H^0(K_S \otimes B^{\otimes n}) \otimes H^0(B) @>\alpha>> H^0(K_S
\otimes B^{\otimes n+1}) \text{ for all }n\geq 2 \cr
H^0(K_S \otimes B^{\otimes m}) \otimes H^0(K_S) @>\beta>>
H^0(K_S^{\otimes 2} \otimes B^{\otimes m}) \text{ for all }m \geq 4 \ .}$$
Let
$C$ be a smooth curve in 
$|K_S|$. The surjectivity of $\beta$ follows by \CM \ns, Lemma 5.2 and the
inequality $B^2 \geq K_S \cdot B$. To see the surjectivity of $\alpha$,
let us set $G= K_S \otimes B^{\otimes n} \otimes \Cal O_C$. Using 
\robs it suffices to show the
surjectivity of 
$$H^0(G) \otimes H^0(B \otimes \Cal O_C)
@>\gamma>> H^0( G \otimes B)  \ .$$ Note that deg$G + \text{ deg} (B \otimes
\Cal O_C) \geq K_S \cdot B + 3B^2 \geq 2(K_S \cdot B + B^2) = 4g(C)-4$.
Since by the inequality $B^2 \geq K_S \cdot B$ and the ampleness of $B$, $h^1(B
\otimes
\Cal O_C) = p_g$ or $p_g-1$, and $p_g \geq 4$, the surjectivity of
$\gamma$ follows from \Pareschilemma \ns.  

To prove the vanishing of $H^{1}(M_{L}^{\otimes
2}\otimes
L^{\prime })$, by  the vanishing
just proved, it suffices to see the surjectivity of 
$$H^0(M_{L}\otimes
L^{\prime }) \otimes H^0(L) \to H^0(M_{L}\otimes
L \otimes L^{\prime }) \ .$$
Using \sfreeobs it is enough to prove the surjectivity of
$$\displaylines{H^0(M_{L}\otimes
L^{\prime }) \otimes H^0(B) \to H^0(M_{L}\otimes
 \otimes L^{\prime } \otimes B) \cr
H^0(M_L \otimes L' \otimes B^{\otimes 2}) \otimes H^0(K_S) \to H^0(M_L \otimes
L' \otimes K_S \otimes B^{\otimes 2}) \ .}$$
The surjectivity of the second family of maps follows by \CM and the same
arguments used for the cohomology vanishing already proven. For the first
family we argue restricting to $C$. By
\robs
\ns, 
\Splemma
\ns, and having in account the already proven surjectivity of
$\gamma$,  we see that it suffices to check the
surjectivity of
$$H^0(M_N \otimes
N^{\prime }) \otimes H^0(B \otimes \Cal O_C) @>\delta >>
H^0(M_{N}\otimes B \otimes N^{\prime }) \ ,$$
where $N = L \otimes \Cal O_C$ and $N'= L' \otimes \Cal O_C$. Now deg$G
\geq 2g(C)-2 + B^2$ and  $B^2 \geq K_S^2 \geq 4$ by Lemma 5.3 and
N\"other's inequality, hence 
$M_G$ is semistable. Then $\delta$ surjects by \Butres \ns.

\vskip .15 truecm

{\it Koszul}: According to \KLaz we need to show that
$M^{h,L}$ is globally generated and that 
$$H^0(M^{h,L})\otimes H^0(L^{\otimes s}) @>\alpha>> H^0(M^{h,L} 
\otimes L^{\otimes s})$$
surjects for all $h \geq 0$ and $s \geq 1$. Let $B'=K_S \otimes B$, ample and
base-point-free by \numfree \ns.  By \sfreeobs \ns, it
suffices to prove the surjectivity of 
$$\displaylines{H^0(M^{h,L}\otimes  B^{\otimes l})\otimes H^0(B) @>\beta_l>>
H^0(M^{h,L} 
\otimes  B^{\otimes l+1})
\text{ for all } l \geq 0 \ , \cr
H^0(M^{h,L}\otimes  B^{\otimes k} \otimes {B'}^{\otimes r})\otimes H^0(B') @>\gamma_r>>
H^0(M^{h,L} 
\otimes  B^{\otimes k} \otimes {B'}^{\otimes r+1})
\text{ for all } k \geq 1, r \geq 0 \ .}$$
We explain in some detail one of the 
border cases:
$$H^0(M^{h,L})\otimes H^0(B) @>\beta>>
H^0(M^{h,L} 
\otimes B)$$
and leave the others to the reader. 
The proof of the surjectivity of
$\beta$ goes by induction and as in \EKoszul, it is
convenient to prove the vanishing of $H^1(M^{h,L}\otimes
B^*)$ at the same time. If
$h=0$, the surjectivity follows from the arguments sketched in the first
part of the proof. Assume the statement to be true for $h-1$.
Consider the sequence 
$$\displaylines{
H^0(M^{h-1,L})\otimes H^0(L \otimes B^*) @>\delta>>
H^0(M^{h-1,L}\otimes L \otimes B^* )\cr
\to H^1(M^{h,L} \otimes B^*) \to H^1(M^{h-1,L})\otimes 
H^0(L \otimes B^*)\ .}$$
The multiplication map $\delta$ is surjective by induction 
hypothesis. The group $H^1(M^{h-1,L})$ vanishes also by induction
hypothesis, therefore $H^1(M^{h,L} \otimes
B^*)=0$. Now since $H^1(\Cal O_X)=0$,  in order to
see the surjectivity of
$\beta$ we proceed as in
\text{\EKoszul \ns,} $B$ playing the role $B_1$ plays there and
restricting to a smooth curve $C \in |B|$, which plays
the same role as $\frak b_1$. In order to obtain the
inequalities needed to apply 
\slopelemma and
\Butres \ns, note that $K_S \otimes B \otimes \Cal O_C=K_C$ and
recall that $B^2 \geq B \cdot K_S$ and deg$(K_S \otimes \Cal O_C) = K_S^2 \geq 4$ by
N\"other's formula. 
$\square$
\vskip .2 truecm
As a corollary we obtain an improvement on another result by 
Ciliberto (cf. \text{\Ciliberto
\ns).} As in Remark 5.7 and \fiveK the hypothesis on the ampleness of
$K_S$ can be removed, and we state the corollary without it:

\proclaim{Corollary 5.15} 
Let $S$ be   a regular 
surface of general type with 
$K_{S}$ 
base-point-free. Let $p_g \geq 4$. 
Let 
$L=K_{S}^{\otimes p+2+l}$. 
Then, if $p \geq 1$,  the image of $S$ by $|L|$ is projectively
normal, its ideal is generated by quadratic equations
and  the homogeneous coordinate ring is
Koszul.
\endproclaim

We end up the section  with  a generalization to higher syzygies 
of Corollary 5.15:

\proclaim{Theorem  5.16} 
Let $S$ be   a regular 
surface of general type with 
$K_{S}$ 
base-point-free. Let $p_g \geq 4$. 
Let 
$L=K_{S}^{\otimes p+2+l}$ and 
$L^{\prime }=K_{S}^{\otimes p+2+k}$. 
Then, if $p \geq 1$,  $H^1(M_{L}\otimes
L^{\prime })=H^{1}(M_{L}^{\otimes
p+1}\otimes
L^{\prime })=0$ for all $k,l\geq 0$. Moreover, if $p
\geq 1$ the image of $S$ by $|L|$ is projectively
normal, its ideal is generated by quadratic equations
and the resolution of the homogeneous coordinate ring is
linear until the $p$th stage.
\endproclaim

\noindent {\it Sketch of proof.}
If $K_S$ is ample, the vanishings
of
$H^1(M_{L}\otimes L^{\prime })$ and $H^1(M_{L}^{\otimes 2}\otimes L^{\prime })$
follows from Theorem 5.14 and if $K_S$ is not ample, they follow by
the same reasoning used for Theorem 5.14, arguing as in Remark 5.7. 
The proof for $p >1$ follows now by induction. We argue as is the  proof  of 
\regulargtNp \ns. We use \sfreeobs \ns, \robs \ns,
\Splemma \ns, and \Butres in similar fashion.

Lastly, the statement about the syzygies of the
resolution of the pluricanonical models follow from
\GLlemma and from the vanishings just proven.
$\square$

\heading References \endheading
\roster

\item
"\Bombieri" E. Bombieri, {Canonical models of surfaces of general type},
Inst. Hautes \'Etudes Sci. Publ. Math., {\bf 42} (1973) 171-219.

\item 
"\Butler" D. Butler, {\it Normal generation of vector bundles over
a curve}, J. Differential 
Geometry {\bf 39} (1994) 1-34.

\item"\Catanese" F. Catanese, {\it Canonical rings and special surfaces of
general type}, Proceedings of the Summer Research Institute on Algebraic
Geometry at Bowdoin 1985, Part 1, AMS 1987, 175-194. 

\item"\Ciliberto" C. Ciliberto, {\it Sul grado dei generatori dell'anello
canonico di una superficie di tipo generale}, Rend. Sem. Mat. Univ. Politecn.
Torino {\bf 41} (1983), 83-111.

\item"\CosDol" F.R. Cossec \& I.V. Dolgachev,  Enriques Surfaces I,
Birkh\"auser, 1989.

\item"\EL" L. Ein \& R. Lazarsfeld, {\it Koszul cohomology and syzygies of
projective
 varieties}, Inv. Math. {\bf  111} (1993), 51-67.

\item"\Busto" G. Fern\'andez del Busto, {\it A Matsusaka-type
theorem on surfaces}, J. of Algebraic Geometry {\bf 5}, (1996),
513-520. 

\item"\FV" M. Finkelberg \& A. Vishik, {\it The coordinate ring of a general
curve of genus $g \geq 5$ is Koszul}, J. of Algebra {\bf 162}, (1993), 535-539.

\item "\GPone" F.J. Gallego \& B.P. Purnaprajna {\it Normal
presentation on elliptic ruled surfaces}, J. of Algebra {\bf 186}, 
(1996),
597-625.

\item"\GPtwo" \hbox{\leaders \hrule  \hskip .6 cm}\hskip .05 cm , {\it Higher
syzygies of elliptic ruled surfaces}, J. of Algebra {\bf 186}, (1996),
626-659.

\item"\GPthree" \hbox{\leaders \hrule  \hskip .6 cm}\hskip .05 cm , 
{\it Vanishing theorems and syzygies for K3 surfaces and Fano
varieties}, preprint.

\item"\GPfour" \hbox{\leaders \hrule  \hskip .6 cm}\hskip .05 cm , second part of
this manuscriprt.
%{\it Very ampleness and higher syzygies for Calabi-Yau threefolds},
%preprint.

\item"\Greentwo" M. Green,
{\it Koszul cohomology and the geometry of projective varieties},
J. Differential Geometry {\bf
19} (1984) 125-171.

\item"\Greenthree" \hbox{\leaders \hrule  \hskip .6 cm}\hskip .05 cm ,
{\it Koszul cohomology and the geometry of projective varieties II},
J. Differential Geometry {\bf
20} (1984) 279-289.

\item"\Greenfour" \hbox{\leaders \hrule  \hskip .6 cm}\hskip .05 cm ,
{\it Koszul cohomology and  geometry}, Lectures on Riemann Surfaces,
 World Scientific
Press, Singapore, 1989, 177-200.

\item"\Hommaone" Y. Homma, {\it Projective normality and the
defining equations of ample 
invertible sheaves on elliptic ruled
surfaces with $e \geq 0$}, Natural Science 
Report, Ochanomizu
Univ. {\bf 31} (1980) 61-73.

\item"\Hommatwo" 
\hbox{\leaders \hrule  \hskip .6 cm}\hskip .05 cm , {\it Projective
normality and the defining equations of an elliptic ruled 
surface with negative
invariant}, Natural Science Report, Ochanomizu
Univ. {\bf 33} (1982) 17-26.

\item"\Kaw" Y. Kawamata, {\it The cone of curves of algebraic
varieties}, Ann. of Math. (2) {\bf 119} (1984), 603-633.

\item"\Kempf" G. Kempf, {\it Projective coordinate rings of Abelian
varieties}, Algebraic Analysis, Geometry and Number Theory, the
Johns Hopkins Univ. Press, 1989, 225-235.

\item"\Lazarsfeld" R. Lazarsfeld, {\it A sampling of vector bundle techniques in
the study of linear series}, Lectures on Riemann Surfaces, World Scientific
Press, Singapore, 1989, 500-559.

\item"\Miyaoka" Y. Miyaoka, {\it The Chern class and Kodaira
dimension of a minimal variety}, 
Algebraic Geometry --Sendai 1985,
Advanced Studies in Pure Math., 
Vol. 10, 
North-Holland,
Amsterdam, 449-476.

\item "\Mumford" D. Mumford
{\it Varieties defined by quadratic equations}, Corso CIME in 
Questions on Algebraic Varieties,
Rome, 1970, 30-100.

\item "\Pareschione" G. Pareschi, {\it Koszul algebras associated to
adjunction bundles}, J. of Algebra 
{\bf 157} (1993) 161-169.

\item"\Pareschitwo" \hbox{\leaders \hrule  \hskip .6 cm}\hskip .05 cm , {\it
Gaussian maps and multiplication maps on certain projective varieties},
Compositio Math. {\bf 98} (1995), 219-268.

\item"\PP" G. Pareschi \& B.P. Purnaprajna, {\it Canonical ring of a 
curve
is Koszul:
A simple proof}, to appear in Illinois J. of Math.

\item"\Reider" I. Reider, {\it Vector bundles of rank $2$  and
linear systems on an algebraic surface}, 
Ann. of Math. (2) {\bf
127} (1988) 309-316.

\item"\Sho" V.V. Shokurov, {\it A nonvanishing theorem}, Math.
USSR-Izv. {\bf 26} (1985) 591-604.

\endroster

\newpage

\hsize= 16.85 truecm
\vsize= 21.75 truecm
\hoffset= -.15 truecm
\voffset= 1 truecm

\def\robss{\GPfour \ns, Observation 2.3 }
\def\Splemmaa{Lemma 2.9 }
\def\ns{\hskip -.15 truecm}
\def\sequence{(2.1) }
\def\Np{Theorem 2.4 }
\def\CM{\Mumford \ns, Theorem 2 }
\def\sfreeobs{Observation 1.2 }
\def\ngfourone{(1.4.1) }
\def\Npreg{Corollary 1.1 }
\def\ngthreeB{Theorem 1.4 }
\def\Greenthm{\Green \ns, Theorem 3.9.3 }
\def\Ksequence{Lemma 3.4 }
\def\conj{Conjecture 1.9 }
\def\robs{Observation 1.3 }
\def\EKoszul{Theorem 3.5 }
\def\Koszul{Theorem 2.7 }
\def\Butres{\Butler \ns, Proposition 2.2 }
\def\foursections{Example 1.5 }
\def\Splemma{Lemma 2.3 }
\def\GLlemma{Theorem 2.2 }
\def\Pareschilemma{\Pareschi \ns, Corollary 4 }

\def\ACGH{[ACGH] }
\def\Butler{[B] }
\def\Ci{[C] }
\def\EL{[EL1] }
\def\ELtwo{[EL2] }
\def\Fujita{[F] }
\def\GPthree{[GP1] }
\def\GPfour{[GP2] }
\def\Green{[G1] }
\def\Greentwo{[G2] }
\def\GL{[GL] }
\def\Kaw{[K] }
\def\Lazarsfeld{[L] }
\def\Mumford{[M] }
\def\OP{[OP] }
\def\Pareschi{[P] }
\def\SD{[S-D] }

\topmatter
\title Part 2: Very Ampleness and higher syzygies for 
Calabi-Yau threefolds
\endtitle
\address{Francisco Javier Gallego: Dpto. de \'Algebra,
 Facultad de Matem\'aticas,
 Universidad Complutense de Madrid, 28040 Madrid,
Spain}\endaddress
\email{gallego\@eucmos.sim.ucm.es}\endemail
\address{ B.P.Purnaprajna:
Dept. of Mathematics,
  University of Missouri, 
Columbia MO 65211}\endaddress
\email{purna\@math.missouri.edu}\endemail
\endtopmatter
\document
\headline={\ifodd\pageno\rightheadline \else\leftheadline\fi}
\def\rightheadline{\tenrm\hfil \eightpoint 
 \hfil\folio}
\def\leftheadline{\tenrm\folio\hfil \eightpoint F.J. GALLEGO \& B.P. PURNAPRAJNA \hfil}

\heading Introduction \endheading

In  this article we prove results on very ampleness, projective normality and
higher syzygies for  Calabi-Yau threefolds.

In the first section  we  prove optimal results
on very ampleness and projective normality for powers of ample and
base-point-free line bundles. Let $X$ be a Calabi-Yau threefold and let $B$ be an
ample and base-point-free line bundle on $X$. The main results of Section 1 can
be summarized in the two following theorems (for a stronger statement of Theorem
2, see Theorem 1.7):

\proclaim{Theorem 1 (cf.  \ngthreeB \ns)}  
The line bundle $B^{\otimes 3}$ is very ample and $|B^{\otimes
3}|$ embeds $X$ as a projectively normal variety  if and only if the morphism
induced by $|B|$ does not map $X \ 2:1$ onto $\bold P^3$.
\endproclaim

\proclaim{Theorem 2}
The line bundle $B^{\otimes 2}$ is very ample and $|B^{\otimes
2}|$ embeds $X$ as a projectively normal variety  if
 $|B|$ does not map $X$ onto a variety of minimal
degree other than $\bold P^3$ nor  maps $X \ 2:1$ onto $\bold P^3$.   
\endproclaim 

A  Calabi-Yau threefold is the three-dimensional version of a K3 surface
and  Theorems 1 and 2 are analogues of  the well known results of
St. Donat (see \SD \ns) for K3 surfaces.  Precisely, for a K3 surface $S$ and an
ample and base-point-free line bundle $B$ on $S$, St. Donat proved the following: 
\roster 
\item
$B^{\otimes 2}$ is very ample and 
$|B^{\otimes 2}|$ embeds $S$ as a projectively normal variety if and only  if the
morphism induced by $|B|$ does not map $S \ 2:1$ onto $\bold P^2$.
\item
$B$ is very ample and 
$|B|$ embeds $S$ as a projectively normal variety if   $|B|$ does not map
$S$ onto a variety of minimal degree nor  maps $X \ 2:1$ onto $\bold P^2$.
\endroster
\vskip .2 truecm
As corollaries of Theorems 1 and 2 and results of Ein, Lazarsfeld, Fujita and
Kawamata  on global generation  on smooth
threefolds, we obtain bounds very closed to Fujita's conjecture.
Precisely  we show the following: 
\proclaim{Corollary 1 (cf. Corollary 1.10)} 
Let $X$ be  a smooth Calabi-Yau threefold
and let $A$ be an ample line bundle. Let $L
=A^{\otimes n}$. 
If $n \geq 8$, then $L$  is very ample and $|L|$ embeds $X$ as a projectively
normal variety. Moreover, if
$A^3 > 1$ and $n \geq 6$, then $L$  is very ample and $|L|$ 
embeds $X$ as a projectively normal variety. 
\endproclaim

We end Section 1 with a result
regarding very ampleness and projective normality on Calabi-Yau fourfolds.

Section 2 is devoted to the computation of Koszul cohomology groups on Calabi-Yau
threefolds. The work of Mark Green in the 80's connected Koszul cohomology with
the study of equations and free resolutions of projective varieties.  From our Koszul
cohomology computations we obtain  results regarding the equations and higher
syzygies associated to  powers of ample and base-point-free line bundles. We also
study the Koszul property for these bundles (see
\Koszul \ns). Regarding equations and higher syzygies we
prove the following 
\proclaim{Theorem 3 (cf. \Np \ns)} Let $X$ be  a  Calabi-Yau threefold and let
$B$ be an ample and base-point-free line bundle on $X$ such that $|B|$ does not
map $X$ onto $\bold P^3$.  If $n \geq p+2$ and $p \geq 1$, then $B^{\otimes n}$
satisfies property
$N_p$. In particular, if $n \geq 3$, the homogeneous ideal associated to the
embedding given by
$|B^{\otimes n}|$ is generated by quadrics.
\endproclaim
The parallelism between K3 surfaces and Calabi-Yau threefolds goes over to higher
syzygies. In fact Theorem 3 is analogous to the following  result
proved by the authors in
\GPfour
\ns:

{\it Let $S$ be a K3 surface and let $B$ be an ample and base-point-free line bundle
on
$S$ such that $|B|$ does not map $S$ onto $\bold P^2$. If $n \geq p+1$ and $p \geq 1$,
then
$B^{\otimes n}$ satisfies property
$N_p$. }

As a corollary of Theorem 3 we obtain bounds for a power of an ample line bundle
to satisfy property
$N_p$. We show precisely the following
\proclaim{Corollary 2 (cf. Corollary 2.8)}
Let $X$ be  a smooth Calabi-Yau threefold
and let $A$ be an ample line bundle. Let $L
=A^{\otimes n}$. 
If $n \geq 4p+8$, then $L$ satisfies property $N_p$ and the coordinate ring of the 
image of the embedding induced by $|L|$ is Koszul. Moreover, if
$A^3 > 1$ and $n \geq 3p+6$, then $L$ satisfies property $N_p$ and the 
coordinate ring of $X$
is Koszul. 
In particular, if $n \geq 12$, or if $n \geq 9$ and $A^3 >1$, then the ideal associated
to the  embedding induced by
$|L|$ is generated by quadratic equations.
\endproclaim

\headline={\ifodd\pageno\rightheadline \else\leftheadline\fi}
\def\rightheadline{\tenrm\hfil \eightpoint VERY AMPLENESS AND HIGHER
SYZYGIES OF CALABI-YAU THREEFOLDS
 \hfil\folio}
\def\leftheadline{\tenrm\folio\hfil \eightpoint F.J. GALLEGO \& B.P. PURNAPRAJNA \hfil}

The article focuses on smooth Calabi-Yau threefolds for  the sake of simplicity.
However the arguments used also go through for Calabi-Yau threefolds with
 terminal singularities and for Calabi-Yau threefolds with 
canonical singularities. In fact Theorems 1, 2 and 3 hold for  Calabi-Yau
threefolds with canonical singularities. From them we recover and strengthen
results by Oguiso and Peternell (see \OP
\ns). The case of singular Calabi-Yau threefolds is dealt with in the appendix at the
end of the article.
\vskip .2 truecm
We  thank
Dale Cutkosky and Mohan Kumar  for their encouragement as well as
for useful discussions. We also thank Sheldon Katz for his encouragement and
discussions regarding  examples of Calabi-Yau threefolds. We are also grateful to
 Vladimir Masek,  who brought to our attention the work
of Oguiso and Peternell, and to Frank Olaf Schreyer.

\vskip .2 truecm

\noindent {\bf Convention.} Throughout this article we
work over an algebraically closed field   
of characteristic
$0$. 

\vskip .2 truecm
\noindent{\bf Definition.} Let $X$ be a projective variety and let $L$
be a very ample line bundle on $X$. We say that $L$ is normally
generated or that $L$ satisfies the property $N_0$, if $|L|$ embeds $X$ as
a projectively normal variety. We say that $L$ is normally presented or
that $L$ satisfies the property $N_1$ if $L$ satisfies property $N_0$
and, in addition, the homogeneous ideal associated to the embedding  of $X$ given
by
$|L|$ is generated by quadratic equations. We say that $L$ satisfies the
property $N_p$ for $p >1$, if $L$ satisfies property $N_1$ and the free
resolution of the homogeneous ideal of $X$  is linear
until the $(p-1)$th-stage.

\heading 1. Very ampleness and projective normality. 
\endheading
A smooth Calabi-Yau threefold $X$ is smooth projective variety of dimension $3$
with trivial canonical bundle and such that $H^1(\Cal O_X)=0$. In this section we 
study when a power of an ample and base-point-free line bundle $A$
 on a Calabi-Yau threefold is very ample and when its complete linear series
embeds $X$ as a projectively normal variety. We recall  the following corollary of
Theorem 1.3 in \GPfour which can be proven using arguments based upon
Castelnuovo-Mumford regularity and Koszul cohomology:

\proclaim{\Npreg (\GPfour \ns, Corollary 1.6)} Let $X$ be a
smooth Calabi-Yau $m$-fold, and
$B$ an ample and base-point-free line bundle on $X$. If $n \geq
p+m$ and $p
\geq 1$  then
$B^{\otimes n}$ satisfies property $N_p$.
\endproclaim

\Npreg tells us in particular
that if $X$ is a Calabi-Yau threefold and $n \geq 4$, then 
$B^{\otimes n}$ satisfies property $N_0$, i.e., is very ample and
embeds the variety as a projectively normal variety.  The main concern of this
section is dealing with the case $n=2$ (Theorem 1.7) and
$n=3$ (\ngthreeB \ns). For that purpose one has to take into account the particular
properties of Calabi-Yau threefolds.  The strategy to follow will be to find suitable
divisors on the threefold and to translate the questions on   surjectivity of
 multiplication maps on the threefold to questions on surjectivity
of multiplication maps on the divisor. These arguments will be fruitfully  repeated,
eventually reaching the situation in which one confronts the question of
surjectivity of multiplications maps on curves. Thus
 results on surjectivity of maps on curves, like \Butres 
and \Pareschilemma,  and on surfaces, like \Greenthm for surfaces of
general  type (see also \Ci \ns),  will be of great interest to us.  Before we proceed
with the statements and proofs of \ngthreeB and Theorem 1.7, we introduce two
auxiliary tools which will be used throughout:

\proclaim {\sfreeobs}
Let $E$, $L_1$ and $L_2$ be 
coherent sheaves on a variety 
$X.$
Consider the multiplication map of global sections $H^0(E) \otimes H^0(L_1
\otimes L_2) @> \psi >> H^0(E \otimes
L_1 \otimes L_2)$ and the maps 
$$\displaylines{H^0(E) \otimes H^0(L_1
) @> \alpha_1 >> H^0(E \otimes L_1) \text { and }
\cr
H^0(E\otimes L_1) \otimes H^0(L_2
) @> \alpha_2 >> H^0(E \otimes L_1\otimes L_2) \ .}$$
If $\alpha_1$ and $\alpha_2$ are surjective
then $\psi$ is also surjective.
\endproclaim

\proclaim{\robs (\robss \ns)}
Let $X$ be a regular variety (i.e, a variety such that 
$H^1(\Cal{O}
_X)=0).$ Let $E$ be a vector bundle on $X$ and let $C$ be a
divisor such that  
$L$
$=\Cal{O}_X\left( C\right) $ is  a globally generated line bundle and 
$H^1(E\otimes L^{-1})=0.$ If the multiplication map 
$$ H^0(E\otimes \Cal{O}_C)\otimes H^0(L\otimes
\Cal{O}_C)\to H^0(E\otimes L\otimes \Cal{O}_C) \text{ surjects},$$   then the
multiplication map 
$$H^0(E)\otimes H^0(L) \to H^0(E\otimes L) \text{ also surjects.}$$
\endproclaim

Now we are ready to state and prove the following  result which
give necessary and sufficient conditions for
$B^{\otimes 3}$ to satisfy property $N_0$:

\proclaim{\ngthreeB}
Let $X$ be a Calabi-Yau threefold and let $B$ be an ample 
and base-point-free line bundle.
Then $B^{\otimes 3}$ is very ample and $|B^{\otimes
3}|$ embeds $X$ as a projectively normal variety except if
$h^0(B)=4$ and the sectional genus of $B$ is 3, in which case
$B^{\otimes 3}$ is not even very ample.
\endproclaim

\noindent {\it Proof.}
{\it Case 1}: $h^0(B) \geq 5$.
It is enough to see that the map
$$
H^{0}(B^{\otimes 3+k})\otimes H^{0}(B^{\otimes 3+l}) \to
H^{0}(B^{\otimes 6+k+l}) 
$$ surjects   for all $ k,l\geq 0$.
By \sfreeobs it is enough 
to prove a stronger statement, namely,   that the map
$$
H^{0}(B^{\otimes 3+l})\otimes H^{0}(B) \to
H^{0}(B^{\otimes 4+l}) 
$$
surjects   for all $ l\geq 0$. Castelnuovo-Mumford regularity arguments will not
work if
$l=0$, so we consider a smooth divisor $S \in |B|$ and the following
commutative diagram:
$$\matrix
H^0(B^{\otimes 3+l}) \otimes H^0(\Cal O_X)&\hookrightarrow
&H^0(B^{\otimes 3+l})
\otimes H^0(B)&\twoheadrightarrow &H^0(B^{\otimes 3+l})
\otimes H^0(B
\otimes 
\Cal O_S)\cr
\downarrow&&
\downarrow&&\downarrow\cr
H^0(B^{\otimes 3+l} )&\hookrightarrow &H^0(B^{\otimes
4+l})&\twoheadrightarrow &H^0(B^{\otimes 4+l} \otimes \Cal
O_S) \ . 
\endmatrix
$$
The map whose surjectivity we wish to show  is 
the middle vertical map. 
The surjectivity of the left hand side vertical map is
obvious.   Note that $B \otimes \Cal O_S = K_S$. Since $H^1(B^{\otimes 2+l})=0$ for
all $l$, checking the surjectivity of the right hand side reduces
to checking the surjectivity of $$
 H^{0}(K_S^{\otimes 3+l})\otimes H^{0}(K_S) @>\alpha>>
H^{0}(K_S^{\otimes 4+l}) \text{ for all } l\geq 0 \ . 
$$
To see that $\alpha$ surjects we consider now a smooth divisor
$C \in |K_S|$. By
\robs
 and Kodaira vanishing, checking the surjectivity of $\alpha$ reduces to checking
the surjectivity of
$$H^0(\theta^{\otimes 3+l}) \otimes H^0(\theta) @>\beta>>
H^0(\theta^{\otimes 4+l})\ ,$$ 
where $\theta=B\otimes \Cal O_C$ is a theta-characteristic. We can now apply either
\Pareschilemma or \Butres  to show the surjectivity of $\beta$. For instance, to
apply  
\Pareschilemma \ns, we need that either $\theta$ or $\theta^{\otimes 3+l}$ be
very ample, that both
$h^0(\theta)$ and $h^0(\theta^{\otimes 3+l})$ be greater than
or equal to $3$ and that deg $\theta^{\otimes 3+l} + \text{ deg }\theta$ be greater
than or equal to both $3g-3$ and
$4g-1-2h^1(\theta)-2h^1(\theta^{\otimes 3+l})-\text{ Cliff}(C)$. The line bundle
$\theta^{\otimes 3+l}$ is very ample because by Clifford's bound $g(C) \geq 5$, and
the required bounds on the number of linearly independent global sections of
$\theta$ and
$\theta^{\otimes 3+l}$ are also satisfied since $h^0(B) \geq 5$. Finally, the last
condition required follows from   deg $\theta^{\otimes 3+l} + \text{ deg }\theta \geq
4g-4$ and
$h^1(\theta) \geq 3$. 
\vskip .1 truecm
\noindent{\it Case 2}: $h^0(B)=4$.
Let $\pi$ be the morphism induced by $|B|$. Let $C$ be as above.
Since
$B
\otimes
\Cal O_C$ has degree
$g(C)-1$ and it is the pullback of $\Cal O_{\bold P^1}(1)$ for a general
$\bold P^1$ in
$\bold P(H^0(B))= \bold P^3$, the degree $n$ of $\pi$ is $g(C)-1$. In
particular,
$g(C) \geq 3$. If $g(C) = 3$, $B^{\otimes 3} \otimes \Cal O_C =
K_C \otimes \theta$, where $\theta$ has degree $2$. Therefore
the restriction of $B^{\otimes 3}$ to $C$ is not very ample. Now we
treat the case  $g(C) \geq 4$. It suffices  to see the surjectivity of
$$
H^{0}(B^{\otimes 3+l})\otimes H^{0}(B^{\otimes 3+k}) \to
H^{0}(B^{\otimes 6+k+l}) \text{ for all } l,k \geq 0 
\ .$$
The key case is $k=l=0$.
If $l \geq 1$ or $k \geq 1$,  the surjectivity of the above map follows from the
arguments displayed below for the case $k=l=0$, or alternatively from
{\text \sfreeobs
\ns,} 
\Mumford \ns, p. 41, Theorem 2 and Kodaira vanishing Theorem. Therefore we
focus our attention on the case
$l=k=0$. It follows from \sfreeobs that it is enough to check the surjectivity of
$$\displaylines{
H^{0}(B^{\otimes 3})\otimes H^{0}(B^{\otimes 2}) @>\alpha>>
H^{0}(B^{\otimes 5})
\cr
H^{0}(B^{\otimes 5})\otimes H^{0}(B) @>\beta>>
H^{0}(B^{\otimes 6})\ .
}$$
The map $\beta$ surjects by \CM and Kodaira vanishing
Theorem.
Note that we cannot use \sfreeobs again in order to prove the
surjectivity of $\alpha$,  because the
map
$H^0(B^{\otimes 3})\otimes H^0(B)
\to H^0(B^{\otimes 4})$ is actually non surjective, for otherwise
the map  
$$
H^{0}(K_C \otimes \theta)\otimes H^{0}(\theta ) \to
H^{0}(K_C^{\otimes 2})$$ 
  would
also surject, but this  is false by 
 base-point-free pencil trick.
Instead  the surjectivity of $\alpha$ will follow from the
surjectivity of $\gamma$ and $\delta$ in the 
diagram
$$\matrix
H^0(B^{\otimes 2}) \otimes  H^0(B^{\otimes 2}) & \hookrightarrow
& H^0(B^{\otimes 2}) \otimes  H^0(B^{\otimes 3})&
\twoheadrightarrow & H^0(B^{\otimes 2}) \otimes 
H^0(K_S^{\otimes 3}) \cr @VV\gamma V @VV\alpha V @VV\delta
V\cr
 H^0(B^{\otimes 4}) & \hookrightarrow &
H^0(B^{\otimes 5})& \twoheadrightarrow &
H^0(K_S^{\otimes 5})  \ ,
\endmatrix
$$
obtained from the sequence
$$ 0 \longrightarrow B^* \longrightarrow \Cal O_X \longrightarrow \Cal O_S
\longrightarrow 0
\eqno
\ngfourone
\ .$$ To see the surjectivity of $\gamma$ we construct yet another two
similar diagrams arising from $\ngfourone$.  
Since
$H^1(B^{\otimes r})=0$ for all $r \geq 0$, checking the surjectivity of
$\gamma$ reduces to seeing the surjectivity of 
$$\displaylines{H^0(K_S^{\otimes 2}) \otimes H^0(K_S^{\otimes 2})
@>\epsilon>> H^0(K_S^{\otimes 4}) \cr
H^0(K_S^{\otimes 2}) \otimes H^0(K_S)
@>\eta>> H^0(K_S^{\otimes 3})\ .}$$
On the other hand in order to see the surjectivity of $\delta$, again
by the vanishing of $H^1(B^{\otimes r})$ it is enough to check 
the surjectivity of 
$$H^0(K_S^{\otimes 3}) \otimes H^0(K_S^{\otimes 2})
@>\varphi>> H^0(K_S^{\otimes 5}) \ .$$
 For the surjectivity of $\epsilon$, $\eta$ and $\varphi$ we
build commutative diagrams like the one above, now upon the
sequence
$$0\to K_S^* \to \Cal O_S \to \Cal O_C \to 0 \ .$$ 
For instance, to see the surjectivity of $\epsilon$ we would write:
$$\matrix
H^0(K_S^{\otimes 2}) \otimes  H^0(K_S) &
\hookrightarrow & H^0(K_S^{\otimes 2}) \otimes 
H^0(K_S^{\otimes 2})&
\twoheadrightarrow & H^0(K_S^{\otimes 2}) \otimes 
H^0(K_C) \cr @VV  \eta V @VV \epsilon V @VV 
V\cr
 H^0(K_S^{\otimes 3}) & \hookrightarrow &
H^0(K_S^{\otimes 4})& \twoheadrightarrow &
H^0(K_C^{\otimes 2})  \ .
\endmatrix
$$
Since $H^1(K_S^{\otimes r})=0$ for all $r \geq 0$, the surjectivity
of $\epsilon$,  $\eta$ and $\varphi$ will follow from the
surjectivity of the maps
$$\displaylines{H^0(K_C) \otimes H^0(K_C) \to
H^0(K_C^{\otimes 2})  \cr
H^0(K_C) \otimes H^0(\theta) \to H^0(K_C 
\otimes \theta) \text{ and } \cr
H^0(K_C) \otimes H^0(K_C \otimes \theta ) \to 
H^0(K_C^{\otimes 2}
\otimes \theta)
\ .}$$
Recall that 
$g(C)
\geq 4$, 
therefore
$C$ cannot be hyperelliptic as  $|B \otimes \Cal O_C|$ is a
base-point-free pencil of degree $g-1$.  Thus the  first map above is surjective by
N\"other's theorem. For the second, recall that $\theta$ is a
theta-characteristic, that it is base-point-free and that 
$h^0(\theta)=2$. Thus the surjectivity
follows from the base-point-free pencil trick.
Finally the third one follows from  \Pareschilemma \ns. 
$\square$

\vskip .2 truecm
We show now by means of an example that there indeed exist ample and
base-point free line bundles  with four linearly independent global sections and
sectional genus $3$:
\vskip .2 truecm
\noindent {\bf \foursections \ns.}
Let $X$ be the double cover of $\bold P^3$ ramified along a 
smooth degree $8$ surface and let $B$ be the pullback of 
$\Cal O_{\bold P^3}(1)$.
The threefold $X$ is Calabi-Yau, $h^0(B)=4$ and the sectional genus of $B$ is $3$.
\vskip .2 truecm
We now want to know when $B^{\otimes 2}$ is normally
generated. In the study we  carry on we will use a theorem by
M. Green. To apply this theorem we will 
require  the image of the morphism induced by $|B|$ not to be
a variety of minimal degree. For that reason it is interesting to
classify the different kinds of varieties of minimal degree which
can appear in our setting and the structure of the Calabi-Yau threefold and of the
morphism induced by $|B|$. This is done in the next proposition.
 
\proclaim{Proposition 1.6} Let $X$ be a smooth Calabi-Yau
threefold, let $\pi$ be the morphism induced by the complete linear series of an
ample and base-point-free line bundle $B$ on $X$ with $h^0(B)=r+1$, and let $n$ be
the degree of
$\pi$. If the image of
$X$ by $\pi$ is a variety $Y$ of minimal degree, then $n \leq \frac{6r} {r-2}$ and
one of the following  occurs:
\roster
\item  $Y= \bold P^3$.
\item $Y$ is  a smooth quadric hypersurface in 
$\bold P^4$.
\item $Y$ is a smooth rational normal scroll of dimension $3$ in
$\bold P^5$. Then the threefold $X$ is fibered over $\bold
P^1$, and  the general fiber is  either a smooth K3 surface, in which case
$n=2,4,6,8$ or $10$, or a smooth Abelian surface, in which case $n=6,8$ or $10$.
\item $Y$ is a smooth rational normal scroll in $\bold P^r$, $r
\geq 6$, the degree $n$ of $\pi$ is $2$ and $X$ is fibered over
$\bold P^1$ with a smooth K3 surface as a general fiber. The 
restriction  of $B$ to the general fiber of $X$ is hyperelliptic, with
sectional genus $2$, and its complete
linear series  maps the fiber 
$2:1$ onto a general fiber of the scroll.
\item $Y$ is a smooth rational normal scroll in $\bold P^r$, $r
\geq 6$, the degree $n$ of $\pi$ is $6$ and $X$ is fibered over
$\bold P^1$ with a smooth Abelian surface as a general fiber.  The  
restriction  of $B$ to the general fiber of $X$ is a $(1,3)$
polarization,  and its complete
linear series  maps the fiber 
$6:1$ onto a general fiber of the scroll.
\item $Y$ is a singular threefold of minimal degree which is either a
cone over a rational normal curve or a cone over a Veronese
surface. 
\endroster
\endproclaim

{\it Proof}. First we prove  the inequality  
$ n \leq \frac {6r} {r-2} \quad (*)$
\ \  holds when if $Y$ is a variety of
minimal degree.  By Riemann-Roch
$r+1 = h^0(B)
\geq
\frac{1}{6} B^3 + 1 = \frac{1}{6} n(r-2)+1$, so we obtain the inequality. 

Now we describe all the possible types of varieties of minimal degree that may
occur. 
The variety
$Y$ should be  either
$\bold P^3$, a smooth quadric hypersurface in $\bold
P^4$, a singular $3$-dimensional
rational normal scroll in $\bold P^4$, a (possibly singular)
$3$-dimensional rational normal scroll in $\bold P^r$, $r \geq 5$
or a cone in
$\bold P^6$ over a  Veronese surface in $\bold P^5$. We see now  that $Y$ cannot be
a   
$3$-dimensional rational normal scroll singular along a single point. 
In that case $Y$ would admit a small resolution and from that it would follow  that
$X$ can also be obtained by performing small contractions on another variety
$\tilde X$, and hence
$X$ would not be smooth.

We complete now the proof of the proposition describing the
cases when $Y$ is a smooth rational normal scroll. In such case, $r
\geq 5$, $Y$ is fibered
over $\bold P^1$ and so is $X$.
Let $\varphi$ the projection from $Y$ to
$\bold P^1$.  Then the general fibers $X @> \varphi \circ \pi >>
\bold P^1$ are, by adjunction,  either smooth K3 surfaces or smooth Abelian
surfaces. Let us denote by $F$ a
general fiber of $\varphi$, and let $G$ be a general  fiber of
$\varphi
\circ
\pi$. We consider the following sequence:
$$ 0 \longrightarrow H^0(B(-G)) \longrightarrow H^0(B)
\longrightarrow H^0(B \otimes \Cal O_G) \longrightarrow
H^1(B(-G)) \longrightarrow 0 \ .$$
If
$r
\geq 6$, $Y=S(a,b,c)$ (i.e., $Y$ is isomorphic to $\bold P(\Cal E)$, where $\Cal E=\Cal
O(a) \oplus \Cal O(b) \oplus \Cal O(c)$, mapped in projective space by $|\Cal O_{\bold
P(\Cal E)}(1)|$), with
$a
\leq b
\leq c$,
$a
\geq 1$, and
$c
\geq 2$. Let $H$ be the restriction of $\Cal O_{\bold P^r}(1)$ to $Y$. Then
$H(-F)$ is big and globally generated, in particular, big and nef,
and $\pi$ being finite, so is $B(-G)$. Thus by Kawamata-Viehweg,
$H^1(B(-G))=0$, which implies that $|B \otimes \Cal O_G|$ maps $G$
onto $\bold P^2$. If $G$ is a smooth K3 surface, then $(G, B
\otimes \Cal O_G)$ is a genus $2$, hyperelliptic polarized K3
surface, $\pi|G$ is $2:1$ and so is $\pi$. If $G$ is a smooth Abelian
surface, then $B \otimes \Cal O_G$ is a $(1,3)$-polarization, hence
$\pi|G$ has degree $6$ and so has $\pi$.

If $r=5$, consider the sequence
$$0 \longrightarrow H^0(B) \longrightarrow H^0(B \otimes G) \longrightarrow
H^0(B \otimes \Cal O_G) \longrightarrow 0 \ .$$
It follows from Riemann-Roch that $h^0(B \otimes \Cal O_G)= h^0(B \otimes
G)-h^0(B)=\frac{1}{2} B^2\cdot G + \frac{1}{12} G \cdot c_2(X)$. Since $B^2 \cdot G =
n$ and
$G
\cdot c_2(X)$ equals $0$ if $G$ is Abelian and $24$ if $G$ is a K3 surface, we finally
obtain that $n = 2h^0(B \otimes \Cal O_G)$ if $G$ is an Abelian surface and $n =
2h^0(B
\otimes \Cal O_G)-4$ if $G$ is a K3 surface. The inequality (*) completes now the
statement made in (3).
$\square$
\vskip .2 truecm
We return our attention to the normal generation of $B^{\otimes
2}$:

\proclaim{Theorem 1.7} Let $X$ be a smooth Calabi-Yau threefold
and let $B$ be an ample and base-point-free line bundle on $X$.
\roster
\item If the image of $X$ by the morphism $\pi$ induced by the complete
linear series of $B$ is not a variety of minimal degree, i.e., is not one of the six
cases in the list of Proposition 1.6, then
$B^{\otimes 2}$ is very ample and embeds $X$ as a projectively
normal variety. 
\item In case 1 of Proposition 1.6, 
 $B^{\otimes 2}$ is very ample and embeds $X$ as a projectively
normal variety if and only if the sectional genus of $B$ is not $3$.
\item If the degree $n$ of $\pi$ equals $2$ (for instance, in case $4$ of Proposition
1.6),
$B^{\otimes 2}$ is not even very ample.
\endroster
\endproclaim

\noindent{\it Proof.}
We prove (1) first. By hypothesis, $h^0(B) \geq 5$ and
 the image of
$X$ by the morphism induced by $|B|$ is not a variety of
minimal degree.  We want to prove that $B^{\otimes 2}$ satisfies property $N_0$.
We prove instead a more general statement, namely, we show 
that the multiplication map 
$$
H^{0}(B^{\otimes l+2})\otimes H^0(B^{\otimes 2}) \longrightarrow
H^{0}(B^{\otimes l+4})
$$
surjects for all $l\geq 0.$ From \sfreeobs it follows that it suffices to have the
surjectivity of 
$$
H^{0}(B^{\otimes l+2})\otimes H^0(B)@> \alpha >>
H^{0}(B^{\otimes l+3})
$$
 for all $l\geq 0$. The crucial  cases are
$l=0,1$. If
$l
\geq 2$, the surjectivity of
$\alpha$  can be obtained from the same arguments used below for $l=0,1$,
or  from
 Kodaira vanishing and \CM \ns.  Case $l =1$ was already dealt with in the proof of 
\ngthreeB \ns. Thus we focus on case $l=0$, i.e., on the surjectivity of
$$H^0(B^{\otimes 2}) \otimes H^0(B) @>\beta >> H^0(B^{\otimes 3})
\ .$$
We  first use
\robs
\ns. Since $H^1(B)=0$ and by adjunction $B \otimes 
\Cal O_S=
K_S$, it is enough to prove the surjectivity of
$$
H^{0}(K_{S})\otimes H^{0}(K_{S}^{\otimes 2})@>{\delta }>>
H^{0}(K_{S}^{\otimes 3})
\ .$$
Since  the image of $S$ under 
the morphism
 defined by $|K_{S}|$ is not a surface of minimal degree,
$h^{0}(K_{S})= h^0(B) -1 \geq 4$, and
$H^{1}(\Cal{O}_{S})=0,$ by \Green \ns, Theorem 3.9.3,  the map
$\delta$ surjects. 
\vskip .2 truecm
We prove now (2). Recall that  $h^0(B)=4$.  We want to show the surjectivity of 
$$
H^{0}(B^{\otimes 2l})\otimes H^0(B^{\otimes 2})@> \alpha >>
H^{0}(B^{\otimes 2l+2}) 
$$
if the sectional genus of $B$ is greater than $3$.
 If $l=1$, the surjectivity of $\alpha$ was shown in the 
proof of
\ngthreeB \ns. If $l \geq 2$, the surjectivity of $\alpha$ follows from the same
arguments or alternatively from  
\text{\sfreeobs \ns,} Kodaira vanishing and \CM \ns. On the other hand, if the
sectional genus of $B$ is
$3$, the morphism induced by $|B|$ is a $2:1$ cover of $\bold
P^3$, hence a general curve $C$ in $B \otimes \Cal O_S$, where
$S$ is a general divisor in $|B|$, is a hyperelliptic curve. Therefore
$B^{\otimes 2} \otimes \Cal O_C= K_C$ is not very ample.

Finally we prove (3). Since now  the morphism induced by $|B|$ is a 
$2:1$ cover of a rational normal scroll and $C$ is again
hyperelliptic, then $B^{\otimes 2}$ cannot be very ample.
$\square$ 

\vskip .2 truecm
 Looking at Theorem 1.7 it can be seen
that the hyperellipticity of $C$ determines in many cases whether
$B^{\otimes 2}$ satisfies the property $N_0$ or not. For instance,
the fact of $C$ being hyperelliptic forces the image of $X$ by $|B|$
to be a variety of minimal degree. We also have this

\proclaim{Corollary 1.8} Let $X$ be a Calabi-Yau threefold and let
$B$ be an ample and base-point-free line bundle on $X$.  If
$h^0(B)=4$ or if $(X,B)$ is of type $2$ (i.e., the morphism
induced by $|B|$ is generically 2:1 onto its image), 
$B^{\otimes 2}$ satisfies the property $N_0$ if and only if there exists $S \in |B|$
and a smooth curve $C \in |B \otimes \Cal O_S|$ which is non-hyperelliptic. 
\endproclaim 
All the above motivates the following
 
\proclaim{\conj} Let $X$ be a Calabi-Yau threefold and let $B$ be an ample and
base-point-free line bundle. Then
$B^{\otimes 2}$ embeds $X$ as a  projectively normal variety if
and only if there is a smooth non-hyperelliptic curve 
$C$ in
$|B
\otimes \Cal O_S|$  some  $S \in |B|$.
\endproclaim

\vskip .2 truecm

Theorems 1.4 and 1.7 combined with results on
global generation of powers of ample line bundles, such as
Ein and Lazarsfeld's (cf. \ELtwo \ns), Fujita's (cf. \Fujita \ns) and
Kawamata's \text{(cf. \Kaw \ns),}  yield the following  

\proclaim{Corollary 1.10} Let $X$ be  a smooth Calabi-Yau threefold
and let $A$ be an ample line bundle. Let $L
=A^{\otimes n}$. 
If $n \geq 8$, then $L$ satisfies property $N_0$. Moreover, if $A^3
> 1$ and $n \geq 6$, then $L$ satisfies property $N_0$. 
\endproclaim

\noindent {\it Proof.} The line bundle  $A^{\otimes m}$ is
base-point-free if $m \geq 4$ (cf. \ELtwo \ns) and, if $A^3  >1$ and
$m \geq 3$, then $A^{\otimes m}$ is base-point-free (cf. \Fujita
\ns). Using
\CM and
\sfreeobs it is  not difficult  to see that
$A^{\otimes n}$ satisfies property $N_0$ if $n \geq 13$ (if $n \geq
10$ when $A^3 >1$). On the other hand, if $n=2l$ with $l \geq 4$ (with
$l \geq 3$ if $A^3>1$),
$A^{\otimes n}$ satisfies property $N_0$ as a consequence of
Theorem 1.7. Indeed, set $B=A^{\otimes l}$. Riemann-Roch implies
that $h^0(B) \geq 12$ ($
\geq 10$, if $A^3>1$). Then by inequality (*), the degree of the map induced by
$|B|$ is less than or equal to $7$ in both cases. By Proposition 1.6 the image
 of $\pi$ is a rational normal scroll $Y$ maybe singular along a line. Let $F$ be a
general $\bold P^2$ among those contained in $Y$ and let $G=\pi^{-1}(F)$. Then 
deg $\pi = B^2\cdot G  \geq 9 A^2 \cdot G \geq 9$, since $A$ is ample, and this is a
contradiction.  
There are therefore only a few cases still to be checked. If  $A^3 >1$, we still have to
deal with  $L=A^{\otimes 7}$ and
$L=A^{\otimes 9}$. The case $m=9$ follows directly from \text{\ngthreeB \ns.} For
the case
$m=7$, we argue as follows. We use
\text{\sfreeobs
\ns.} The only multiplication maps whose surjectivity cannot be
checked using \CM are 
$$\displaylines{H^0(A^{\otimes 7}) \otimes H^0(A^{\otimes 3})
@>\alpha>> H^0(A^{\otimes 10}) \cr
H^0(A^{\otimes 10}) \otimes H^0(A^{\otimes 4}) @>\beta>>
H^0(A^{\otimes 14}) \ .}$$
To see the surjectivity of $\alpha$ and $\beta$ one can argue as in
Theorem 1.7, reducing the problem eventually to checking the
surjectivity of multiplication maps on suitable curves. The
surjectivity of these maps follows by \Pareschilemma or \Butres \ns.
Finally, if we are in the case when no conditions are
imposed on
$A^3$, then we still have  to deal with $m=9, 11$ and $12$. The
argument is the same as before. 
 $\square$

\vskip .2 truecm
To end this section we prove a result regarding very ampleness and projective
normality on Calabi-Yau fourfolds. Recall that  \Npreg tells among other things that
if
$X$ is a smooth Calabi-Yau fourfold and $B$ is an ample and base-point-free line
bundle, then 
$B^{\otimes n}$ satisfies property
$N_p$ if $n \geq p+4$ and $p \geq 1$. Therefore if $n \geq 5$, $B^{\otimes n}$
satisfies property
$N_1$, and in particular, $B^{\otimes n}$ is very ample and
$|B^{\otimes n}|$ embeds
$X$ as a projectively normal variety. In the following theorem we prove that the
above holds for $B^{\otimes 4}$ under certain   conditions on $B$.

\proclaim{Theorem 1.11} Let $X$ be a smooth Calabi-Yau fourfold and let $B$ be an
ample and base-point-free line bundle such that the morphism induced by $|B|$ is
birational onto the image and $h^0(B) \geq 7$. Then $B^{\otimes 4}$ is very ample
and
$|B^{\otimes 4}|$ embeds $X$ as a projectively normal variety.
\endproclaim

\noindent{\it Proof}.
From \sfreeobs it follows that it suffices to prove the surjectivity of 
$$H^0(B^{\otimes n}) \otimes H^0(B) \longrightarrow H^0(B^{\otimes n+1})$$
for all $n \geq 4$. When $n \geq 5$, this follows from \CM and Kodaira vanishing
theorem. If
$n=4$, we argue like in the proofs of Theorems 1.4 and 1.7. We consider a smooth
curve
$C$ obtained by iteratively taking  hyperplane sections in
$|B|$. Then we use \robs and since $B^{\otimes 4} \otimes \Cal O_C= K_C$, 
the problem is eventually reduced  to checking the
surjectivity of the following map on $C$,
$$H^0(K_C) \otimes H^0(L) @>\alpha >> H^0(K_C \otimes L) \ ,$$
where $L= B \otimes \Cal O_C$.
The line bundle $L$  is ample, base-point-free,
$|L|$ induces a birational morphism from $C$ onto its image, and $h^0(L) \geq
4$, thus the surjectivity of $\alpha$ follows from 
 a theorem of Castelnuovo (cf. \ACGH \ns, page 151) which states that the
map
$S^n H^0(L)
\otimes H^0(K_C) \longrightarrow H^0(K_C \otimes L^{\otimes n})$ surjects for all
$n \geq 0$ under the  conditions satisfied by $L$.

\heading 2. Normal presentation, Koszul rings and higher
syzygies. \endheading

The purpose of this section is to compute 
Koszul cohomology groups on Calabi-Yau threefolds and to apply this computation
to the study of the ring, equations and free resolution of those threefolds. The
connection between Koszul cohomology and syzygies was developed by Green (see
\Greentwo \ns; for a particularly gentle introduction to the subject see
also \Lazarsfeld \ns).   We present now the statement we need for our
purposes. Let
$X$ be a projective variety, and  let $F$ be a globally generated
vector bundle on $X$. We define the bundle $M_F$ as follows:
$$ 0 \to M_F \to H^0(F) \otimes \Cal O_X \to F \to 0 \ .\eqno
\sequence
$$
If $L$ is an ample line bundle on $X$ such that all its positive powers are
nonspecial there exists  the following characterization of  property
$N_p$:

\proclaim {\GLlemma } Let $L$ be an ample, globally generated line
bundle on a variety $X$. If  the cohomology group
\text{$H ^1(\bigwedge ^{p'+1} M_L \otimes
L^{\otimes s})$} vanishes for all $0 \leq p' \leq p$ and all $s \geq
1$,  then
$L$ satisfies the
 property
$N_p$. If in addition  $H^1(L ^{\otimes r}) = 0$, for all $r \geq 1$,
then the above  is a necessary and sufficient condition
for $L$ to satisfy property $N_p$. 
\endproclaim

\vskip .2 truecm
We use this characterization to prove our results on syzygies. For the proof of it we
refer to
\text{\EL
\ns,} Section 1. Since we are working over an
algebraically closed field of characteristic $0$,  for our proofs of higher syzygies
results we will check 
 the vanishings of  
\text{$H ^1( M_L^{\otimes p'+1} \otimes
L^{\otimes s})$} rather than see directly the vanishings of 
\text{$H ^1(\bigwedge ^{p'+1} M_L \otimes
L^{\otimes s})$}. 
\vskip .2 truecm
Before we state the main theorem of this section we state the following lemma
(for the proof, see \GPfour \ns, \Splemmaa \ns):

\proclaim{\Splemma} Let $X$ be a projective variety, let $q$ be
a nonnegative integer and let $F$ be a base-point-free line bundle
on $X$. Let $Q$ be an effective line
bundle on $X$ and let $\frak q$ be a reduced and irreducible member
of
$|Q|$. 
Let $R$ be a line bundle and $G$ a sheaf on $X$ such that
\vskip .1 cm
\item{1.} $\text H^1(F \otimes Q^*)=0$,
  \item {2.}  
$
\text H
^0(M_{(F\otimes \Cal O_\frak q)}^{\otimes q'}
  \otimes
R \otimes \Cal O_\frak q)
\otimes \text H^0(G) \to \text H ^0(M_{F\otimes \Cal O_\frak
q)}^{\otimes q'}
\otimes R \otimes G \otimes \Cal O_\frak q) \ \text {surjects for all}
\  0 \leq q' \leq q$.

\noindent Then, for all $0 \leq q'' \leq q$ and all $0 \leq k' \leq q''$, 
$$\text H^0(M_{F}^{\otimes k}
\otimes M_{(F \otimes \Cal O_\frak q)}^{\otimes q''-k}  \otimes R
\otimes \Cal O_\frak q) \otimes \text H^0(G) \to 
\text
H^0(M_{F}^{\otimes k}
\otimes M_{(F \otimes \Cal O_\frak q)}^{\otimes q''-k} \otimes G \otimes R
\otimes \Cal O_\frak q)  \text{ surjects.}$$ 
\endproclaim 
 
We are now ready the state the following
\proclaim{\Np}
Let $X$ be a Calabi-Yau threefold. Let $B$ be an ample and 
base-point-free divisor with $h^{0}(B)\geq 5.$ Let
$L=B^{\otimes p+2+k}$ and $L'=B^{\otimes p+2+l}$. If $k,l\geq
0$ and $p \geq 1$, 
then $H^{1}(M_{L}^{\otimes p+1}\otimes
L^{\prime })=0$ and $L$ satisfies property $N_p$. 
\endproclaim

\noindent{\it Proof.}
The proof is by induction on $p$. The most important step is $p=1$.  Consider the
sequence
$$\displaylines {H^0(M_L \otimes L') \otimes H^0(L) @> \alpha >> 
H^0(M_L
\otimes L'
\otimes L) \cr 
\to H^1(M_L^{\otimes 2} \otimes L') \to H^1(M_L
\otimes L') \otimes H^0(L) \ .} $$
The last term of the sequence vanishes by \ngthreeB \ns, so
it suffices  to prove the surjectivity of $\alpha$. For that we use
\sfreeobs \ns. We see therefore that it is enough to show that
$$ H^0(M_L \otimes L') \otimes H^0(B) @> \beta >> 
H^0(M_L
\otimes L'
\otimes B)$$ surjects. 
Let
$S$ be a smooth divisor in
$|B|$. The cohomology group
$H^1(M_L
\otimes B)$ vanishes because the map
$\alpha$ surjects, as shown in the proof of \ngthreeB \ns. Thus by \robs it is
enough to show the surjectivity of 
$$H^0(M_L \otimes B \otimes \Cal O_S) \otimes H^0(B \otimes \Cal
O_S)  @> \gamma >> H^0(M_L \otimes B^{\otimes 2} \otimes \Cal
O_S)
\ .$$ 
Applying now   
\Splemma
\ns, we conclude that it suffices to see the surjectivity of 
$$H^0(M_{K_S^{\otimes m}} \otimes K_S^{\otimes n}) \otimes H^0(K_S) @> \delta >>
H^0(M_{K_S^{\otimes m}}
 \otimes K_S^{\otimes n+1})\ , $$  
for all $m,n \geq 3$.  Let $C$ be  a smooth curve in $|K_S|$ and set $G=K_S^{\otimes
m} \otimes \Cal O_C$ and $G'=K_S^{\otimes n} \otimes \Cal O_C$.  We apply 
\robs
and  
\Splemma
\ns. To apply \Splemma we need  to see that 
$$\displaylines {H^0(
G^{\prime }) \otimes H^0(\theta) \longrightarrow
H^0(G' \otimes \theta) \text{ and } \cr
H^0(M_{G} \otimes
G') \otimes H^0(\theta) \longrightarrow
H^0(M_{G} \otimes
G' \otimes \theta) }$$
surject, where $\theta = B \otimes \Cal O_C$ is a theta-characteristic. To see the
surjectivity of the first map, note that deg
$(G' \otimes \Cal O_C)+
\text{ deg }\theta \geq 4g(C)-4$. Since $h^0(B) \geq 5$, then $h^1(\theta) \geq 3$,
so the surjectivity follows by \Butres or \Pareschilemma \ns. To see that the second
map surjects,  note that $K_S^2 \geq 4$ by
N\"other's inequality, and therefore,  deg $G\geq 3g(C)-3
\geq 2g(C)+2$. Thus
$M_{G} \otimes G'$ is semistable by \Butler \ns, Theorem 1.12. We see now that the
slope of
$M_{G}
\otimes G'$ is bigger than $2g(C)$. Since $H^1(G)=0$, 
$$\mu(M_{G})= \frac {-\text{ deg }G} {\text{ deg }G-g(C)} \ ,$$
therefore $$\mu(M_{G} \otimes G') \geq \frac {-\text{ deg }G } {\text{ deg }G-g(C)}
+3g(C)-3 \geq \frac {-\text{ deg }G} {\text{ deg }G-g(C)} + 2g(C)+2 \ ,$$
and the last term of the above sequence of inequalities is bigger than or equal to
$2g(C)+1$. On the other hand
$$\mu(M_{G} \otimes  G')= \frac {-\text{ deg }G} {\text{ deg }G-g(C)} +3g(C)-3 >
3g(C)-5\geq 2g(C)+2g(C)-\text{ deg }(\theta) -2h^1(\theta) \ .$$
Thus the desired surjectivity follows from
\Butres
\ns.

 For
$p >1$, we write a similar sequence: 
$$\displaylines {H^0(M_L^{\otimes p} \otimes L') \otimes H^0(L) @>
\epsilon >>  H^0(M_L
\otimes L'
\otimes L) \cr 
\to H^1(M_L^{\otimes 2} \otimes L') \to H^1(M_L
\otimes L') \otimes H^0(L) \ .} $$
The group $H^1(M_L^{\otimes p}
\otimes L')$  vanishes by induction hypothesis. By \sfreeobs
we only need to show that 
$$ H^0(M_L^{\otimes p} \otimes L') \otimes H^0(B) @>
\eta >>  H^0(M_L
\otimes L'
\otimes B)$$ surjects.  This follows arguing similarly as in the
proof of the surjectivity of $\beta$, using \text{\robs \ns,} \Splemma to reduce the
problem to checking the surjectivity of multiplication maps on $S \in |B|$ first and to
checking the surjectivity of multiplication maps on $C \in |K_S|$ eventually, and
once we are arguing on $C$, the result follows from \Butres \ns. We can
 argue alternatively by induction to see that  
$\eta$ surjects. Indeed, applying \CM  \ns, the vanishing of
$H^1(M_L^{\otimes p}\otimes L' \otimes B^*)$  follows by
induction and the other two vanishings required  follow from
chasing cohomology sequences arising from \sequence and again
using induction. 

 Finally, the fact that
$L$ satisfies property
$N_p$ follows from the vanishings just proven, \ngthreeB and
\GLlemma \ns.
$\square$
\vskip .2 truecm
\Np says in particular that $B^{\otimes n}$
satisfies property $N_1$, i.e., that the image of the embedding induced by
$|B^{\otimes n}|$ is  ideal-theoretically cut out by quadrics, if $h^0(B) \geq 5$ and $n
\geq 3$. The bound imposed on $h^0(B)$ is sharp, since \foursections provides
an example in which $h^0(B)=4$ and $B^{\otimes 3}$ does 
not even satisfy property
$N_0$  (cf.
\ngthreeB
\ns). We present now an example in which $B^{\otimes 3}$ satisfies property
$N_0$, but not property $N_1$:

\proclaim{Example 2.5} Let $X$ be a cyclic triple cover of $\bold P^3$
ramified along a smooth sextic surface and let $B$ be the pullback of $\Cal
O_{\bold P^3}(1)$ to $X$. The threefold $X$ is a Calabi-Yau threefold and
$h^0(B)=4$. By Theorem 1.7, $B^{\otimes 3}$ satisfies property $N_0$. However,
$B^{\otimes 3}$ does not satisfy property $N_1$. 
\endproclaim 

\noindent {\it Proof.} We sketch the proof of the last claim. Assume $L=B^{\otimes
2}$ satisfies
$N_1$. By
\GLlemma the assumption implies 
$$ H^1(\bigwedge ^2 M_L
\otimes L^{\otimes n})=0 \eqno (2.5.1) $$ 
for all $n \geq 1$. Let $S \in |B|$ and let $C$ be a
smooth curve in $|B \otimes \Cal O_C|$. Using
\sequence it can be seen that both 
$H^2(M_L^{\otimes 2} \otimes L^{\otimes n} \otimes B^*)$ and $H^2(M_L^{\otimes 2}
\otimes L^{\otimes n} \otimes B^*\otimes \Cal O_S)$ vanish. Those vanishings
together with (2.5.1) imply that $H^1(\bigwedge ^2 M_L
\otimes L^{\otimes n} \otimes \Cal O_C)=0$. On the other hand there is an
epimorphism between  the vector bundles $M_L \otimes \Cal O_C$ and $M_{L \otimes
\Cal O_C}$ on $C$. Therefore we have  $$ H^1(\bigwedge ^2 M_{L \otimes \Cal
O_C}
\otimes L^{\otimes n} )=0 \eqno (2.5.2)$$ for all $n \geq 1$. It is a well known result
by Castelnuovo that a line bundle of degree greater than or equal to $2g+1$ on a
smooth curve satisfies property $N_0$. The curve
$C$ has genus
$4$ and $L \otimes \Cal O_C$ has degree  $9$, hence $L
\otimes \Cal O_C$ satisfies $N_0$. Thus it would follow from (2.5.2) that $L \otimes
\Cal O_C$ satisfies also property $N_1$. But this contradicts  a result by Green and
Lazarsfeld (cf. \GL \ns), which says that a line bundle cannot satisfy property $N_1$
if it is the tensor product of the canonical bundle on $C$ and an effective line bundle
of degree
$3$, as is the case of $L \otimes \Cal O_C$. Therefore the original assumption
(2.5.1) is false and $L$ does not satisfies property $N_1$. 
$\square$

\vskip .2 truecm
\noindent{\bf Remark 2.6.} In \Np we dealt with the vanishings needed for
property $N_p$, but in fact the arguments used yield  more general cohomology
vanishings: 

\noindent {\it Let $X$ be a smooth Calabi-Yau threefold and let $B$ be an ample
and base-point-free line bundle such that $h^0(B)\geq 5$. Then
$H^1(M_{B^{\otimes n_1}} \otimes \cdots \otimes M_{B^{\otimes n_{p+1}}} \otimes
B^{\otimes n})=0$ for all $n \geq p+2$, $n_1 \geq 3$ and $n_2, \dots n_{p+1} \geq 1$}.

\vskip .2 truecm
We show now that the line bundle $L$ of \Np
embeds the Calabi-Yau threefold as a variety with a Koszul
coordinate ring. 

\proclaim{\Koszul}
Let $X$ be a Calabi-Yau threefold. Let $B$ be an ample and 
base-point-free divisor with $h^{0}(B)\geq 5.$ Let
$L=B^{\otimes p+2+k}$ and $L'=B^{\otimes p+2+l}$. If $k,l\geq
0$ and $p \geq 1$, 
then  the coordinate
ring of the  image of the embedding induced by $|L|$ is Koszul. 
\endproclaim

\noindent{\it Sketch of proof.}
We follow the same philosophy used in  other
proofs in this article. The claim follows from results
regarding the Koszul property proven for surfaces of general
type in
\GPfour \ns. Precisely it follows as a corollary of \GPfour \ns, Theorem 5.14,
 using
\GPfour \ns, \Ksequence \ns,  and \robs with the same strategy
used to prove \GPfour \ns, \EKoszul \ns. $\square$ 

\vskip .2 truecm
As we did in Section 1, we obtain the following corollary for
powers of ample line bundles. The proof is analogous to the proof
of Corollary 1.10.
\proclaim{Corollary 2.8}
Let $X$ be  a smooth Calabi-Yau threefold
and let $A$ be an ample line bundle. Let $L
=A^{\otimes n}$. 
\roster
\item If $n \geq 4p+8$, then $L$ satisfies property $N_p$. Moreover, if
$A^3 > 1$ and $n \geq 3p+6$, then $L$ satisfies property $N_p$. 
\item If $n \geq 12$ or if $A^3 > 1$ and $n \geq 9$, the coordinate ring of the  image
of the embedding induced by
$|A^{\otimes n}|$ is Koszul.
\endroster
\endproclaim

\heading Appendix: Singular Calabi-Yau threefolds.
\endheading

Throughout the previous part of this article we have been
concerned only with smooth Calabi-Yau threefolds for
reasons of simplicity. In this appendix we show that our
arguments can be adapted without much difficulty to
 canonical Calabi-Yau threefolds and that  our main theorems hold indeed
for them. 

There were only two instances in the proofs of Theorem 1.4 and 1.7 when the
assumption of the nonsingularity of $X$ was used. The first
of them  was when we wanted to guarantee the vanishing of
$H^1(B^{\otimes n})$ for an ample line bundle $B$ and all
$n \geq 0$. This vanishing holds as well for   Calabi-Yau
threefolds with   canonical singularities.   
The second  was to find, firstly a smooth surface $S$ in $|B|$, and
secondly a smooth curve $C$ in $|B \otimes \Cal O_S|$. 
If $X$ has canonical  singularities, it is not possible
in general to find a smooth surface $S$ in $|B|$, since $S$ could have
(at worst) rational double points, but it is possible to find a smooth
curve $C$ in $|B \otimes \Cal O_S|$ using a Bertini-type argument.
Thus the only troublesome point is the use of Green's theorem.
However this result can  still be applied to  $S$ general in $|B|$ if $X$ has
canonical singularities, since we can apply it to the
desingularization $\tilde S$ of $S$, for $S$ and $\tilde S$  have the same canonical
ring. The upshot of all this is that 
\text{{\bf \ngthreeB} \ns,} {\bf Theorem 1.7}, {\it and} {\bf \Np}  {\it hold for
canonical  Calabi-Yau threefolds}, having only in account in the case of 
Theorem 1.7 that there is another case to add to  Proposition 1.6 (6), namely
the image of
$X$ being a cone over a smooth $2$-dimensional rational normal scroll.

As in the end of Sections 1 and 2,  we state now corollaries regarding powers of
ample line bundles. We use for this purpose a generalization of Ein and
Lazarsfeld's result on base-point-freeness, carried out by Oguiso and Peternell.
They prove among other things (cf.
\OP
\ns, Theorems I(2) and II(2)) that $A^{\otimes n}$ is base-point-free if
$n \geq 5$ and $A$ is an ample line bundle on a 
 Calabi-Yau threefold with $\bold Q$-factorial terminal and if $n \geq 7$ and $A$ is
an ample line bundle on a canonical  Calabi-Yau threefold.
As corollaries we recover their results on normal generation of
powers of ample line bundles (see \OP \ns, Theorems I(3) and
II(3)) and generalize them to normal presentation and higher
syzygies. We point out that \OP \ns,  Theorem 3 can also be
recovered as corollary of our Theorem 1.7.

\proclaim{Corollary A.1 (\OP \ns, Theorem I, (3))}
Let $X$ be   Calabi-Yau threefold with $\bold Q$-factorial terminal 
singularities and let
$A$ be an ample line bundle on $X$. If $n \geq 10$, then  
$A^{\otimes n}$ satisfies property $N_0$.
\endproclaim

\proclaim{Corollary A.2 (\OP \ns, Theorem II,  (3))}
Let $X$ be   Calabi-Yau threefold with canonical 
singularities and let
$A$ be an ample line bundle on $X$. If $n \geq 14$, then  
$A^{\otimes n}$ satisfies property $N_0$.
\endproclaim

As corollaries of \Np \ns, we obtain:

\proclaim{Corollary A.3}
Let $X$ be   Calabi-Yau threefold with $\bold Q$-factorial terminal 

singularities and let
$A$ be an ample line bundle on $X$. If $n \geq 5p+10$, then  
$A^{\otimes n}$ satisfies property $N_p$. Furthermore, if $p \geq 1$, the
coordinate ring of the  image of the embedding induced by $|A^{\otimes n}|$ is
Koszul.
\endproclaim

\proclaim{Corollary A.4}
Let $X$ be   Calabi-Yau threefold with canonical 

singularities and let
$A$ be an ample line bundle on $X$. If $n \geq 7p+14$, then  
$A^{\otimes n}$ satisfies property $N_p$. Furthermore, if $p \geq 1$, the
coordinate ring of the  image of the embedding induced by $|A^{\otimes n}|$ is
Koszul.
\endproclaim

\heading References \endheading
\roster

\item"\ACGH" E. Arbarello, M. Cornalba, P.A. Griffiths \& J. Harris,  Geometry of
Algebraic Curves, Volume I, Springer-Verlag 1985.

\item 
"\Butler" D. Butler, {\it Normal generation of vector bundles over
a curve}, J. Differential 
Geometry {\bf 39} (1994) 1-34.

\item"\Ci" C. Ciliberto, {\it Sul grado dei generatori dell'anello
canonico di una superficie di tipo generale}, Rend. Sem. Mat. Univ. Politecn.
Torino {\bf 41} (1983), 83-111.

\item"\EL" L. Ein \& R. Lazarsfeld, {\it Koszul cohomology and syzygies of
projective
 varieties}, Invent. Math. {\bf  111} (1993), 51-67.

\item"\ELtwo" \hbox{\leaders \hrule  \hskip .6 cm}\hskip .05 cm , {\it Global
generation of pluricanonical and adjoint linear series on smooth projective
threefolds}, J.  Amer. Math. Soc. {\bf 6} (1993), 875-903.

\item"\Fujita" T. Fujita, {\it Remarks on Ein-Lazarsfeld criterion of spannedness of
adjoint bundles on polarized threefolds}, preprint.

\item"\GPthree" F.J. Gallego \& B.P. Purnaprajna, 
{\it Vanishing theorems and syzygies for K3 surfaces and Fano
varieties}, preprint.

\item"\GPfour" \hbox{\leaders \hrule  \hskip .6 cm}\hskip .05 cm , first part of this
manuscript.
%{\it Projective normality and syzygies of algebraic surfaces},
%preprint.

\item"\Green" M. Green, {\it The canonical ring of a variety of general type}, 
Duke Math. J. {\bf 49},(1982), 1087-1113.

\item"\Greentwo" \hbox{\leaders \hrule  \hskip .6 cm}\hskip .05 cm ,
{\it Koszul cohomology and the geometry of projective varieties},
J. Differential Geometry {\bf
19} (1984) 125-171.

\item"\GL" M. Green \& R. Lazarsfeld, {\it Some results on the
syzygies of finite sets and  
algebraic curves}, Compositio
Math. {\bf 67} (1989) 301-314.

\item"\Kaw" Y. Kawamata, {\it On Fujita's freeness conjecture for $3$-folds and
$4$-folds}, Math. Ann. {\bf 308}, (1997), 491-505.

\item"\Lazarsfeld" R. Lazarsfeld, {\it A sampling of vector bundle techniques in
the study of linear series}, Lectures on Riemann Surfaces, World Scientific
Press, Singapore, 1989, 500-559.

\item "\Mumford" D. Mumford
{\it Varieties defined by quadratic equations}, Corso CIME in 
Questions on Algebraic Varieties,
Rome, 1970, 30-100.

\item"\OP" K. Oguiso \& T. Peternell, {\it On polarized canonical Calabi-Yau
threefolds}, Math. Ann. {\bf 301}, (1995), 237-248. 

\item"\Pareschi" G. Pareschi, {\it
Gaussian maps and multiplication maps on certain projective varieties},
Compositio Math. {\bf 98} (1995), 219-268.

\item"\SD" Saint-Donat, {\it On projective models of K3 surfaces}, Amer. J. Math. {\bf
96} (1974), 602-639.

\endroster

\enddocument